\documentclass[12pt]{amsart}
\usepackage{verbatim}
\usepackage{amsmath}
\usepackage{amssymb}
\usepackage{amsfonts}
\usepackage{amscd}
\usepackage{graphics}
\usepackage{pictex}
\usepackage{color}

\numberwithin{equation}{section}
 
\newtheorem{theorem}[equation]{Theorem}
\newtheorem{lemma}[equation]{Lemma}
\newtheorem{conjecture}[equation]{Conjecture}

\newtheorem{corollary}[equation]{Corollary}

\newtheorem*{theorem*}{Theorem}
 
\theoremstyle{definition}

\newtheorem{example}[equation]{Example}
 
\theoremstyle{remark}
\newtheorem{remark}[equation]{Remark}

\newcommand{\la}{{\langle}}
\newcommand{\ra}{{\rangle}}
\renewcommand{\O}{{\mathcal O}}
\newcommand{\X}{{\mathfrak X}}

\newcommand{\J}{{\mathcal J}}

\newcommand{\C}{{\mathcal C}}
\newcommand{\D}{{\mathcal D}}
\newcommand{\F}{{\mathcal F}}

\newcommand{\M}{{\mathcal M}}

\newcommand{\tK}{{\tilde{K}}}

\newcommand{\Dbar}{{\bar{D}}}
\newcommand{\Pbar}{{\bar{P}}}

\newcommand{\onto}{\twoheadrightarrow}  

\def\CC{{\mathbb C}}
\def\QQ{{\mathbb Q}}

\def\FF{{\mathbb F}}
\def\NN{{\mathbb N}}
\def\RR{{\mathbb R}}
\def\ZZ{{\mathbb Z}}
\def\PP{{\mathbb P}}

\DeclareFontFamily{OT1}{rsfs}{}
\DeclareFontShape{OT1}{rsfs}{n}{it}{<-> rsfs10}{}
\DeclareMathAlphabet{\mathscr}{OT1}{rsfs}{n}{it}

\newenvironment{notation}[0]{%
  \begin{list}%
    {}%
    {\setlength{\itemindent}{0pt}
     \setlength{\labelwidth}{4\parindent}
     \setlength{\labelsep}{\parindent}
     \setlength{\leftmargin}{5\parindent}
     \setlength{\itemsep}{0pt}
     }%
   }%
  {\end{list}}

\newcommand{\Div}{\operatorname{Div}}

\newcommand{\Jac}{\operatorname{Jac}}
\newcommand{\Prin}{\operatorname{Prin}}

\newcommand{\Ker}{\operatorname{Ker}}
\newcommand{\Image}{\operatorname{Im}}

\newcommand{\Gal}{\operatorname{Gal}}
\newcommand{\Kbar}{{\bar K}}

\newcommand{\Pic}{\operatorname{Pic}}
\newcommand{\Spec}{\operatorname{Spec}}

\newcommand{\cl}{{\rm cl}}
\newcommand{\an}{{\rm an}}
\newcommand{\nr}{{\rm nr}}

\begin{document}

\title[Specialization of Linear Systems]{Specialization of Linear Systems from Curves to Graphs}

\date{July 5, 2007}


\subjclass[2000]{Primary 14H25, 05C99; Secondary 14H51, 14H55.}

\author{Matthew Baker}
\address{School of Mathematics\\
Georgia Institute of Technology\\
Atlanta, Georgia 30332-0160\\
USA}
\email{mbaker@math.gatech.edu}

\begin{abstract}
We investigate the interplay between linear systems on curves
and graphs in the context of specialization of divisors on an
arithmetic surface.  We also provide some applications of our results 
to graph theory, arithmetic geometry, and tropical geometry.  
\end{abstract}

\thanks{
I am very grateful to Marc Coppens for pointing out an error in an earlier version of
this paper; his comments led to numerous improvements in the present version.
I would also like to thank Serguei Norine for many enlightening conversations,
my student Adam Tart for performing some very useful computer calculations, and Michael Kerber
for stimulating discussions concerning tropical curves.
I am also grateful to Kevin Buzzard, Bas Edixhoven, Xander Faber, and Ezra Miller
for their comments on earlier versions of this paper.
Finally, I would like to express my gratitude to Brian Conrad for 
contributing Appendix~\ref{DeformationTheorySection}.  \\
\indent{T}he author's work was supported in part by NSF grant DMS-0600027, and
Brian Conrad's work was partially supported by NSF grant DMS-0600919.}
                                       
\maketitle


\section{Introduction}
\label{IntroSection}

\subsection{Notation and terminology}
\label{section:Notation}

We set the following notation, which will be used throughout the paper unless otherwise
noted.

\begin{notation}
\item[$G$] 
a graph, by which we will mean a finite, unweighted, connected multigraph without loop edges.
We let $V(G)$ (resp. $E(G)$) denote the set of vertices (resp. edges) of $G$.
\item[$\Gamma$]
a metric graph (see \S\ref{MetricGraphSection} for the definition).
\item[$\Gamma_{\QQ}$]
the set of ``rational points'' of $\Gamma$ (see \S\ref{MetricGraphSection}).
\item[$R$]
a complete discrete valuation ring with field of fractions $K$ and algebraically
closed residue field $k$.
\item[$\Kbar$]
a fixed algebraic closure of $K$.
\item[$X$]
a smooth, proper, geometrically connected curve over $K$.
\item[$\X$]
a proper model for $X$ over $R$.
For simplicity, we assume unless otherwise stated that $\X$ is regular, that 
the irreducible components of $\X_k$ are all
smooth, and that all singularities of $\X_k$ are ordinary double points.
\end{notation}

\medskip

Unless otherwise specified, by a {\em smooth curve} we will always mean a
smooth, proper, geometrically connected curve over a field, and by
an {\em arithmetic surface} we will always mean a proper flat scheme of dimension 2 over a
discrete valuation ring whose generic fiber is a smooth curve.  We will usually, but not always,
be working with {\em regular} arithmetic surfaces.  
A {\em model} for a smooth curve $X/K$ is an arithmetic surface $\X/R$ whose generic fiber is $X$.
An arithmetic surface $\X$ is called {\em semistable} if its special fiber $\X_k$ is reduced
and has only ordinary double points as singularities.
If in addition the irreducible components of $\X_k$ are all smooth 
(so that there are no components with self-crossings), we will 
say that $\X$ is {\em strongly semistable}.

\subsection{Overview}

In this paper, we show that there is a close connection between
linear systems on a curve $X/K$ and linear systems, in the sense 
of \cite{BakerNorine}, on the dual graph $G$ of a regular 
semistable model $\X/R$ for $X$.
A brief outline of the paper is
as follows.  In \S\ref{SpecializationSection}, we prove
a basic inequality -- the ``Specialization Lemma'' -- which says
that the dimension of a linear system can only go up
under specialization from curves to graphs (see Lemma~\ref{SpecializationLemma} for the precise statement) .  
In the rest of the paper, we explore various
applications of this fact, illustrating along the way a fruitful interaction between 
divisors on graphs and curves.
The interplay works in both directions: for example, 
in \S\ref{ApplicationSection} we use Brill-Noether theory for curves to prove and/or conjecture 
some new results about graphs (c.f.~Theorem~\ref{GraphBrillNoetherTheorem} and Conjecture~\ref{BrillNoetherConj})
and, in the other direction, in \S\ref{WeierstrassSection} we use the
notion of Weierstrass points on graphs to gain new insight into
Weierstrass points on curves (c.f.~Corollaries~\ref{WeierstrassSpecializationCor} and \ref{WeierstrassBananaCor}).

\medskip

Another fruitful interaction which 
emerges from our approach is a ``machine'' for transporting certain theorems about curves 
from classical to tropical\footnote{See \cite{MikhalkinICM} for an introduction to tropical geometry.} algebraic geometry.
The connection goes through the theory of arithmetic surfaces, by way of the deformation-theoretic result proved
in Appendix~\ref{DeformationTheorySection}, and uses the approximation method introduced in \cite{GathmannKerber} 
to pass from $\QQ$-graphs to arbitrary metric graphs, and finally to tropical curves.
As an illustration of this machine, we prove an analogue for tropical curves 
(c.f.~Theorem~\ref{TropicalBrillNoetherTheorem} below)
of the classical fact that
if $g,r,d$ are nonnegative integers for which the Brill-Noether number
$\rho(g,r,d) := g - (r+1)(g-d+r)$ is nonnegative, then on 
every smooth curve $X/\CC$ of genus $g$ there is a divisor $D$ with $\dim |D| = r$ and $\deg(D) \leq d$.
We also prove, just as in the classical case of algebraic curves, 
that there exist Weierstrass points on every tropical curve of genus $g \geq 2$
(c.f.~Theorem~\ref{WeierstrassExistenceTheorem}).

\medskip

We conclude the paper with two appendices which can be read independently of the rest of the paper.
In Appendix~\ref{RaynaudAppendix}, 
we provide a reformulation of certain parts of 
Raynaud's theory of ``specialization of the Picard functor'' \cite{RaynaudPicard}
in terms of linear systems on graphs.  We also point out some useful consequences of Raynaud's results
for which we do not know any references.
Although we do not actually use Raynaud's results in the body of the paper, 
it should be useful for future work on the interplay between curves and graphs
to highlight the compatibility between Raynaud's theory and our notion of linear equivalence on graphs.
Appendix~\ref{DeformationTheorySection}, written by Brian Conrad,
discusses a result from the deformation theory of stable marked curves 
which implies that every finite graph occurs as the dual graph of a regular semistable model 
for some smooth curve $X/K$ with totally degenerate special fiber.
This result, which seems known to the experts but for which we could not find a suitable reference, 
is used several times throughout the main body of the paper.

\subsection{Divisors and linear systems on graphs}
\label{GraphDivisorSection}

By a {\em graph}, we will always mean a finite, connected multigraph without 
loop edges.  

Let $G$ be a graph, and 
let $V(G)$ (resp. $E(G)$) denote the set of vertices (resp. edges) of $G$.
We let $\Div(G)$ denote the free abelian group on $V(G)$, and refer to elements of $\Div(G)$ as
{\em divisors} on $G$.  
We can write each divisor on $G$ as
$D = \sum_{v \in V(G)} a_v (v)$ with $a_v \in \ZZ$,
and we will say that $D \geq 0$ if $a_v \geq 0$ for all $v \in V(G)$.
We define $\deg(D) = \sum_{v \in V(G)} a_v$ to be the {\em degree} of $D$.
We let $\Div_+(G) = \{ E \in \Div(G) \; : \; E \geq 0 \}$ 
denote the set of {\em effective} divisors on $G$, and we let
$\Div^0(G)$ denote the set of divisors of degree zero on $G$.
Finally, we let $\Div_+^d(G)$ denote the set of effective divisors of degree $d$ on $G$.

Let $\M(G)$ be the group of $\ZZ$-valued functions on $V(G)$, and
define the Laplacian operator
$\Delta : \M(G) \to \Div^0(G)$ by
\[
\Delta(\varphi) = \sum_{v \in V(G)} \sum_{e = vw \in E(G)} 
\left(\varphi(v) - \varphi(w) \right)(v).
\]
We let
\[
\Prin(G) = \Delta(\M(G)) \subseteq \Div^0(G)
\]
be the subgroup of $\Div^0(G)$ consisting of {\em principal divisors}.

\medskip

If $D,D' \in \Div(G)$, we write $D \sim D'$ if $D - D' \in \Prin(G)$,
and set
\[
|D| = \{ E \in \Div(G) \; : \; E \geq 0 \textrm{ and } E \sim D \}.
\]
We refer to $|D|$ as the {\em (complete) linear system} associated
to $D$, and call divisors $D,D'$ with $D \sim D'$ {\em linearly equivalent}.

\medskip

Given a divisor $D$ on $G$, define
$r(D) = -1$ if $|D| = \emptyset$, and otherwise set
\[
r(D) = \max \{ k \in \ZZ \; : \; |D - E| \neq \emptyset \; \forall \; E \in \Div_+^k(G) 
\}.
\]
Note that $r(D)$ depends only on the linear equivalence class of $D$,
and therefore can be thought of as an invariant of the complete
linear system $|D|$.

When we wish to emphasize the underlying graph $G$, we will
sometimes write $r_G(D)$ instead of $r(D)$.

\medskip

We define the {\em canonical divisor} on $G$ to be 
\[
K_G = \sum_{v \in V(G)} (\deg(v) - 2)(v).
\]
We have $\deg(K_G) = 2g-2$, where $g = |E(G)| - |V(G)| + 1$ 
is the {\em genus} (or {\em cyclomatic number}) of $G$.

\medskip

The Riemann-Roch theorem for graphs, proved in \cite[Theorem 1.12]{BakerNorine}, is the following assertion:

\begin{theorem}[Riemann-Roch for graphs]
\label{GraphRiemannRoch}
Let $D$ be a divisor on a graph $G$.  Then
\[
r(D) - r(K_G - D) = \deg(D) + 1 - g.
\]
\end{theorem}

One of the main ingredients in the proof of Theorem~\ref{GraphRiemannRoch} 
is the following result \cite[Theorem 3.3]{BakerNorine}, which is also quite useful 
for computing $r(D)$ in specific examples.
For each linear ordering $<$ of the vertices of $G$, we define a corresponding divisor $\nu \in \Div(G)$
of degree $g-1$ by the formula
\[
\nu = \sum_{v \in V(G)}(|\{e = vw \in E(G) \; : \; w < v \}| - 1)(v).
\]

Then:

\begin{theorem}
\label{OrderTheorem}
For every $D \in \Div(G)$, exactly one of the following holds:
\begin{itemize}
\item[(1)] $r(D) \ge 0$; or
\item[(2)] $r(\nu - D) \ge 0$ for some divisor $\nu$ associated to a linear ordering $<$ of $V(G)$.
\end{itemize}
\end{theorem}

\subsection{Subdivisions, metric graphs, and $\QQ$-graphs}
\label{MetricGraphSection}

By a {\em weighted graph}, we will mean a graph in which each edge is
assigned a positive real number called the {\em length} of the edge.
Following the terminology of \cite{BF}, a 
{\em metric graph} (or {\em metrized graph}) is 
a compact, connected metric space $\Gamma$ which arises by viewing 
the edges of a weighted graph $G$ as line segments.
Somewhat more formally, a metric graph should be thought of as corresponding to an {\em equivalence class} of 
weighted graphs, where two weighted graphs $G$ and $G'$ are {\em equivalent} if they admit a common refinement.
(A {\em refinement} of $G$ is any weighted graph obtained by subdividing the edges of $G$ in a length-preserving fashion.)
A weighted graph $G$ in the equivalence class corresponding to $\Gamma$ is called a {\em model} for $\Gamma$.
Under the correspondence between equivalence classes of weighted graphs and metric graphs,
after choosing an orientation, 
each edge $e$ in the model $G$ can be identified with the real interval $[0,\ell(e)] \subseteq \Gamma$.

\medskip

We let $\Div(\Gamma)$ denote the free abelian group on the points of the
metric space $\Gamma$, and refer to elements of $\Div(\Gamma)$ as
{\em divisors} on $\Gamma$.  
We can write an element $D \in \Div(\Gamma)$ as 
$D = \sum_{P \in \Gamma} a_P (P)$ with $a_P \in \ZZ$ for all $P$ and $a_P = 0$
for all but finitely many $P$.
We will say that $D \geq 0$ if $a_P \geq 0$ for all $P \in \Gamma$.
We let $\deg(D) = \sum_{P \in \Gamma} a_P$ be the {\em degree} of $D$,
we let $\Div_+(\Gamma) = \{ E \in \Div(\Gamma) \; : \; E \geq 0 \}$ 
denote the set of {\em effective} divisors on $\Gamma$, and we let
$\Div^0(\Gamma)$ denote the subgroup of divisors of degree zero on $\Gamma$.
Finally, we let $\Div_+^d(\Gamma)$ denote the set of effective divisors of degree $d$ on $\Gamma$.

\medskip

Following \cite{GathmannKerber}, a {\em $\QQ$-graph} 
is a metric graph $\Gamma$ having a model $G$ 
whose edge lengths are rational numbers.
We call such a model a {\em rational model} for $\Gamma$.
An ordinary unweighted graph $G$ can be thought of as a $\QQ$-graph whose edge lengths are all 1.
We denote by $\Gamma_{\QQ}$ the set of points of $\Gamma$ whose distance from 
every vertex of $G$ is rational; we call elements of $\Gamma_{\QQ}$ {\em rational points}
of $\Gamma$.  It is immediate that the set $\Gamma_{\QQ}$ does not depend on the choice 
of a rational model $G$ for $\Gamma$.
We let $\Div_{\QQ}(\Gamma)$ be the free abelian group on $\Gamma_{\QQ}$, and refer
to elements of $\Div_{\QQ}(\Gamma)$ as {\em $\QQ$-divisors} on $\Gamma$.

\medskip

A {\em rational function} on a metric graph $\Gamma$ is a continuous, piecewise affine function $f : \Gamma \to \RR$,
all of whose slopes are integers.
We let $\M(\Gamma)$ denote the space of rational functions on $\Gamma$.
The {\em divisor} of a rational function $f \in \M(\Gamma)$ is defined\footnote{Here we follow 
the sign conventions from \cite{BF}.  
In \cite{GathmannKerber}, the divisor of $f$ is defined to be the negative of the one we define here.}
to be
\[
(f) = - \sum_{P \in \Gamma} \sigma_P(f) (P),
\]
where $\sigma_P(f)$ is the sum of the slopes of $\Gamma$ in all directions emanating from $P$.
We let $\Prin(\Gamma) = \{ (f) \; : \; f \in \M(\Gamma) \}$ be the subgroup of $\Div(\Gamma)$ consisting
of {\em principal divisors}.  It follows from Corollary 1 in \cite{BF} that
$(f)$ has degree zero for all $f \in \M(\Gamma)$, i.e., $\Prin(\Gamma) \subseteq \Div^0(\Gamma)$.

If $\Gamma$ is a $\QQ$-graph, we denote by $\Prin_{\QQ}(\Gamma)$ the group of principal divisors supported on $\Gamma_{\QQ}$.

\begin{remark}
\label{LaplacianRemark}
As explained in \cite{BF}, if we identify a rational function $f \in \M(\Gamma)$ with
its restriction to the vertices of any model $G$ for which $f$ is affine along each edge of $G$, then 
$(f)$ can be naturally identified with the {\em weighted combinatorial Laplacian} $\Delta(f)$ of $f$ on $G$.
\end{remark}

If $D,D' \in \Div(\Gamma)$, we write $D \sim D'$ if $D - D' \in \Prin(\Gamma)$,
and set
\[
|D|_{\QQ} = \{ E \in \Div_{\QQ}(\Gamma) \; : \; E \geq 0 \textrm{ and } E \sim D \}
\]
and
\[
|D| = \{ E \in \Div(\Gamma) \; : \; E \geq 0 \textrm{ and } E \sim D \}.
\]

It is straightforward using Remark~\ref{LaplacianRemark} to show that if $G$ is a graph and $\Gamma$ is the
corresponding $\QQ$-graph all of whose edge lengths are 1, then two divisors $D,D' \in \Div(G)$ are equivalent on $G$
(in the sense of \S\ref{GraphDivisorSection}) if and only they are equivalent on $\Gamma$ in the sense just defined.

\medskip

Given a $\QQ$-graph $\Gamma$ and a $\QQ$-divisor $D$ on $\Gamma$, define
$r_{\QQ}(D) = -1$ if $|D|_{\QQ} = \emptyset$, and otherwise set
\[
r_{\QQ}(D) = \max \{ k \in \ZZ \; : \; |D - E|_{\QQ} \neq \emptyset \; \forall \; E \in \Div_{\QQ}(\Gamma) \textrm{ with } E \geq 0, \, 
\deg(E) = k \}.
\]
Similarly, given an arbitrary metric graph $\Gamma$ and a divisor $D$ on $\Gamma$, we define
$r_{\Gamma}(D) = -1$ if $|D| = \emptyset$, and otherwise set
\[
r_{\Gamma}(D) = \max \{ k \in \ZZ \; : \; |D - E| \neq \emptyset \; \forall \; E \in \Div_+^k(\Gamma) \}.
\]

\medskip

Let $k$ be a positive integer, and let $\sigma_k(G)$ be the graph
obtained from the (ordinary unweighted) graph $G$ by subdividing each 
edge of $G$ into $k$ 
edges.  
A divisor $D$ on $G$ can also be thought of as a divisor on 
$\sigma_k(G)$ for all $k\geq 1$ in the obvious way.
The following recent combinatorial result \cite{HKN}, 
which had been conjectured by the author, relates the quantities $r(D)$
on $G$ and $\sigma_k(G)$:

\begin{theorem}[Hladky-Kr{\'a}l-Norine]
\label{SubdivisionTheorem}
Let $G$ be a graph.  
If $D \in \Div(G)$, then for every integer $k \geq 1$, we have
\[
r_G(D) = r_{\sigma_k(G)}(D).
\]
\end{theorem}

When working with metric graphs, if $\Gamma$ is the metric graph corresponding
to $G$ (in which every edge of $G$ has length 1), then we will usually think of each edge of 
$\sigma_k(G)$ as having length $1/k$; in this way, each finite graph 
$\sigma_k(G)$ can be viewed as a model for the same underlying metric $\QQ$-graph $\Gamma$.

It is evident from the definitions that $r_{\QQ}(D)$ and $r_\Gamma(D)$ do not change if the length of
every edge in (some model for) $\Gamma$ is multiplied by a positive integer $k$.
Using this observation, together with Proposition 2.4 of \cite{GathmannKerber},
one deduces from Theorem~\ref{SubdivisionTheorem} the following result:

\begin{corollary}
\label{SubdivisionCorollary}
If $G$ is a graph and $\Gamma$ is the corresponding metric $\QQ$-graph
in which every edge of $G$ has length 1, then for every divisor
$D$ on $G$ we have
\begin{equation}
\label{eq:RankInvariance}
r_G(D) = r_{\QQ}(D) = r_{\Gamma}(D).
\end{equation}
\end{corollary}

By Corollary~\ref{SubdivisionCorollary}, we may unambiguously write
$r(D)$ to refer to any of the three quantities appearing in
(\ref{eq:RankInvariance}).

\begin{remark}
Our only use of Theorem~\ref{SubdivisionTheorem} in this paper, other than the notational convenience
of not having to worry about the distinction between $r_{\QQ}(D)$ and $r_G(D)$, 
will be in Remark~\ref{WeierstrassSubdivisionRemark}. 
In practice, however, Theorem~\ref{SubdivisionTheorem} is 
quite useful, since without it there is no obvious way to calculate the quantity 
$r_{\QQ}(D)$ for a divisor $D$ on a metric $\QQ$-graph $\Gamma$.

On the other hand, we {\em will} make use of the equality $r_{\QQ}(D) = r_{\Gamma}(D)$ 
given by \cite[Proposition 2.4]{GathmannKerber} when we develop a machine
for deducing theorems about metric graphs from the corresponding results for $\QQ$-graphs
(c.f.~\S\ref{TropicalBrillNoetherSection}).
\end{remark}

Finally, we recall the statement of the Riemann-Roch theorem for
metric graphs.
Define the {\em canonical divisor} on $\Gamma$ to be 
\[
K_\Gamma = \sum_{v \in V(G)} (\deg(v) - 2)(v)
\]
for any model $G$ of $\Gamma$.  It is easy to see that $K_\Gamma$ is independent of the
choice of a model $G$, and that $\deg(K_\Gamma) = 2g-2$, where $g = |E(G)| - |V(G)| + 1$ 
is the {\em genus} (or {\em cyclomatic number}) of $\Gamma$.

\medskip

The following result is proved in 
\cite{GathmannKerber} and \cite{MikhalkinZharkov}:

\begin{theorem}[Riemann-Roch for metric graphs]
\label{MetricGraphRiemannRoch}
Let $D$ be a divisor on a metric graph $\Gamma$.  Then
\begin{equation}
\label{eq:MetricRR}
r_{\Gamma}(D) - r_{\Gamma}(K_\Gamma - D) = \deg(D) + 1 - g.
\end{equation}
\end{theorem}

By Corollary~\ref{SubdivisionCorollary}, there is a natural
`compatibility' between Theorems~\ref{GraphRiemannRoch}
and \ref{MetricGraphRiemannRoch}.

\subsection{Tropical curves}
\label{TropicalCurveSection}

Tropical geometry is a relatively recent and highly active area of research, and in dimension one it is closely
connected with the theory of metric graphs as discussed in the previous section.
For the sake of brevity, we adopt a rather minimalist view of tropical curves in this paper -- the
interested reader should see \cite{GathmannKerber,MikhalkinICM,MikhalkinZharkov} for motivation and a more extensive discussion.

\medskip

Following \cite[\S{1}]{GathmannKerber}, we define a {\em tropical curve}
to be a ``metric graph with possibly unbounded ends''.  In other words, the only difference between a tropical curve $\tilde{\Gamma}$ and a metric
graph $\Gamma$ is that we allow finitely many edges of $\tilde{\Gamma}$ to have infinite length.\footnote{Unlike some other definitions 
in the literature, the definition of a tropical curve from \cite{GathmannKerber} allows vertices of valence 1 and 2, and
requires that there is a ``point at infinity'' at the end of each unbounded edge.}
We assume that each unbounded edge of $\tilde{\Gamma}$ can
be identified with the extended real interval $[0,\infty]$ in such a way that the $\infty$ end of the edge has valence 1.
Note that every metric graph is also a tropical curve.

\medskip

One can define divisors, rational functions, and linear equivalence for tropical curves exactly as we have done in \S\ref{MetricGraphSection}
for metric graphs; see \cite[\S1 and Definition~3.2]{GathmannKerber} for details.  
(The only real difference is that one must allow a rational function to
take the values $\pm\infty$ at the unbounded ends of $\tilde{\Gamma}$.)
Using the tropical notion of linear equivalence, one defines $r_{\tilde{\Gamma}}(D)$ for divisors on tropical curves 
just as we defined $r_{\Gamma}(D)$ in \S\ref{MetricGraphSection} for divisors on metric graphs.
With these definitions in place, the Riemann-Roch formula (\ref{eq:MetricRR}) holds in the context of tropical curves,
a result which can be deduced easily from Theorem~\ref{MetricGraphRiemannRoch} (see \cite[\S{3}]{GathmannKerber} for details).

\section{The Specialization Lemma}
\label{SpecializationSection}

In this section, we investigate the behavior of the quantity $r(D)$ under
specialization from curves to graphs.  In order to make this precise,
we first need to introduce some notation and background facts
concerning divisors on arithmetic surfaces.

\subsection{The specialization map}

Let $R$ be a complete discrete valuation ring with field of fractions $K$ and algebraically closed
residue field $k$.
Let $X$ be a smooth curve over $K$,
and let $\X / R$ be a strongly semistable regular model for $X$ with special fiber $\X_k$
(see \S\ref{section:Notation}).  
We let $\C = \{ C_1,\ldots,C_n \}$ be the set of irreducible components of $\X_k$.

\medskip

Let $G$ be the {\em dual graph} of $\X_k$, i.e.,
$G$ is the finite graph whose vertices $v_i$ correspond to the irreducible components $C_i$ of $\X_k$,
and whose edges correspond to intersections between these components (so that there is one edge
between $v_i$ and $v_j$ for each point of intersection between $C_i$ and $C_j$).
The assumption that $\X$ is strongly semistable implies that $G$ is well-defined and has
no loop edges.

\medskip

We let $\Div(X)$ (resp. $\Div(\X)$) be the group of Cartier divisors on $X$ (resp. on $\X$); 
since $X$ is smooth and $\X$ is regular, Cartier divisors on $X$ (resp. $\X$) are the same
as Weil divisors.  Recall that when $K$ is perfect, $\Div(X)$ can be identified with the group
of $\Gal(\Kbar/K)$-invariant elements of the free abelian group $\Div(X(\Kbar))$ on $X(\Kbar)$.

We also let $\Prin(X)$ (resp. $\Prin(\X)$) denote the group of principal Cartier divisors on $X$ (resp. $\X$).

\medskip

There is a natural inclusion $\C \subset \Div(\X)$, and 
an intersection pairing 
\[
\begin{aligned}
\C \times \Div(\X) &\to \ZZ \\
(C_i, \D) &\mapsto (C_i \cdot \D),
\end{aligned}
\]
where $(C_i \cdot \D) = \deg(\O_{\X}(\D)|_{C_i})$.

The intersection pairing gives rise to a homomorphism 
$\rho : \Div(\X) \to \Div(G)$
by the formula
\[
\rho(\D) = \sum_{v_i \in G} (C_i \cdot \D) (v_i). 
\]
(If we wish to emphasize the dependence on the ground field $K$, we will sometimes write
$\rho_K(\D)$ instead of $\rho(\D)$.)
We call the homomorphism $\rho$ the {\em specialization map}.
By intersection theory, the group $\Prin(\X)$ is contained in the kernel of $\rho$.

\medskip

The Zariski closure in $\X$ of an effective divisor on $X$ is a Cartier divisor.  Extending by
linearity, we can associate to each $D \in \Div(X)$ a Cartier divisor $\cl(D)$ on $\X$, which we
refer to as the {\em Zariski closure} of $D$.
By abuse of terminology, we will also denote by $\rho$ the
composition of $\rho : \Div(\X) \to \Div(G)$ with the map ${\rm cl}$.
By construction, if $D \in \Div(X)$ is {\em effective}, then $\rho(D)$ is an effective
divisor on $G$.
Furthermore, $\rho: \Div(X) \to \Div(G)$ is a degree-preserving homomorphism.

\medskip

A divisor $\D \in \Div(\X)$ is called {\em vertical} if it is supported on $\X_k$, and {\em horizontal}
if it is the Zariski closure of a divisor on $X$.
If $\D$ is a vertical divisor, 
then $\rho(\D) \in \Prin(G)$; this follows from the fact that
$\rho(C_i)$ is the negative Laplacian of the characteristic function of the vertex $v_i$.
Since every divisor $\D \in \Div(\X)$ can be
written uniquely as $\D_h + \D_v$ with $\D_h$ horizontal 
and $\D_v$ vertical, it follows that 
$\rho(\D)$ and $\rho(\D_h)$ are linearly equivalent divisors on $G$.

Consequently, if $D \in \Prin(X)$, then although the horizontal divisor $\D := \cl(D)$ 
may not belong to $\Prin(\X)$, it differs from a principal divisor $\D' \in \Prin(\X)$ 
by a vertical divisor $\F \in \Div(\X)$ for which
$\rho(\F) \in \Prin(G)$.  Thus we deduce the following basic fact:

\begin{lemma}
\label{PrincipalSpecializationLemma}
If $D \in \Prin(X)$, then $\rho(D) \in \Prin(G)$.
\end{lemma}

\medskip

When $D$ corresponds to a Weil divisor $\sum_{P \in X(K)} n_P(P)$ supported on $X(K)$, there is another,
more concrete, description of the map $\rho$.
Since $\X$ is regular, each point $P \in X(K)=\X(R)$ specializes to a nonsingular point of $\X_k$, and hence to a 
well-defined irreducible component $c(P) \in \C$, which
we may identify with a vertex $v(P) \in V(G)$.
Then by \cite[Corollary~9.1.32]{Liu}, we have
\begin{equation}
\label{ConcreteRho}
\rho(D) = \sum_P n_P (v(P)).
\end{equation}

\begin{remark}
\label{SurjectivityRemark}
Since the natural map from $X(K)=\X(R)$ to the smooth locus of $\X_k(k)$ is surjective (see, e.g.,
Proposition 10.1.40(b) of \cite{Liu}), it follows from
(\ref{ConcreteRho}) that $\rho : \Div(X) \to \Div(G)$ is surjective.  In fact,
this implies the stronger fact that the restriction of $\rho$ to $X(K) \subseteq \Div(X)$ is surjective.
\end{remark}

\subsection{Behavior of $r(D)$ under specialization}
\label{SpecializationSubsection}

Let $D \in \Div(X)$, and let $\Dbar = \rho(D) \in \Div(G)$ be its specialization to $G$.  
We want to compare the dimension of the complete linear system $|D|$ on $X$ 
(in the sense of classical algebraic geometry) with the quantity $r(\Dbar)$ defined in 
\S\ref{GraphDivisorSection}.
In order to do this, we first need some simple facts about linear systems on curves.  

\medskip

We temporarily suspend our convention that $K$ denotes the field of fractions of a discrete
valuation ring $R$, and allow $K$ to be an arbitrary field, with $X$ still denoting
a smooth curve over $K$.
We let $\Div(X(K))$ denote the free abelian group on $X(K)$, which
we can view in a natural way as a subgroup of $\Div(X)$.
For $D \in \Div(X)$, let $|D| = \{ E \in \Div(X) \; : \; E \geq 0, \; E \sim D \}$.
Set $r(D) = -1$ if $|D| = \emptyset$, and otherwise put
\[
r_{X(K)}(D) := \max \{ k \in \ZZ \; : \; |D - E| \neq \emptyset \; \forall \; E \in \Div_+^k(X(K)) \}.\\
\]

\begin{lemma}
\label{CurveRankLemma}
Let $X$ be a smooth
curve over a field $K$, and assume that $X(K)$ is infinite.  
Then for $D \in \Div(X)$, we have
$r_{X(K)}(D) = \dim_K L(D) - 1$, where $L(D) = \{ f \in K(X) \; : \; (f) + D \geq 0 \} \cup \{ 0 \}$.
\end{lemma}

\begin{proof}
It is well-known that $\dim L(D-P) \geq \dim L(D) - 1$ for all $P \in X(K)$.
If $\dim L(D) \geq k+1$, it follows that for any points $P_1,\ldots,P_k \in X(K)$
we have $\dim L(D-P_1-\cdots-P_k) \geq 1$, so that 
$L(D-P_1-\cdots-P_k) \neq (0)$ and $|D-P_1-\cdots-P_k| \neq \emptyset$.

Conversely, we prove by induction on $k$ that if $\dim L(D) = k$, then
there exist $P_1,\ldots,P_k \in X(K)$ such that $L(D-P_1-\cdots-P_k) = (0)$,
i.e., $|D-P_1-\cdots-P_k| = \emptyset$.
This is clearly true for the base case $k = 0$.  
Suppose $\dim L(D) = k \geq 1$, and choose a nonzero rational 
function $f \in L(D)$.
Since $f$ has only finitely many zeros and $X(K)$ is infinite, there exists
$P = P_1 \in X(K)$ for which $f(P) \neq 0$.  It follows that 
$L(D-P) \subsetneq L(D)$, so that $\dim L(D-P) = k-1$.  By induction,
there exist $P_2,\ldots,P_{k} \in X(K)$ such that 
$|D-P-P_2\cdots-P_k| = \emptyset$, which proves what we want.
\end{proof}

Since $\dim L(D)$ remains constant under base change by an arbitrary field extension
$K'/K$, we conclude:
 
\begin{corollary}
\label{CurveRankCor}
Let $X$ be a smooth
curve over a field $K$, and assume that $X(K)$ is infinite.  
Let $K'$ be any extension field.  Then for $D \in \Div(X)$, we have
$r_{X(K)}(D) = r_{X(K')}(D)$.
\end{corollary}

In view of Lemma~\ref{CurveRankLemma} and Corollary~\ref{CurveRankCor}, we will simply write $r_X(D)$, or even 
just $r(D)$, to denote the quantity $r_{X(K)}(D) = \dim_K L(D) - 1$.

\begin{lemma}
\label{CurveGenericPointLemma}
Let $X$ be a smooth
curve over a field $K$, and assume that $X(K)$ is infinite.  
If $D \in \Div(X)$, then $r(D - P) \geq r(D) - 1$ for all $P \in X(K)$, 
and if $r(D) \geq 0$, then $r(D-P) = r(D) - 1$ for some $P \in X(K)$.
\end{lemma}

\begin{proof}
Let $k = r(D)$.
The result is clear for $r(D) \leq 0$,
so we may assume that $k \geq 1$.  If $P = P_1,P_2,\ldots,P_{k} \in X(K)$ are arbitrary, then since 
$r(D) \geq k$, we have
\[
|D - P - P_2 - \cdots - P_k| \neq \emptyset,
\]
and therefore $r(D - P) \geq k-1$.
Also, since $r(D) = k$, it follows that there exist $P=P_1,P_2,\ldots,P_{k+1} \in X(K)$ such that
\[
|D - P - P_2 - \cdots - P_{k+1}| = \emptyset,
\]
and therefore $r(D-P) \leq k-1$ for this particular choice of $P$.
\end{proof}

The same proof shows that an analogous result holds in the context of graphs:

\begin{lemma}
\label{GraphGenericPointLemma}
Let $G$ be a graph, and let $D \in \Div(G)$.
Then $r(D - P) \geq r(D) - 1$ for all $P \in V(G)$, 
and if $r(D) \geq 0$, then $r(D-P) = r(D) - 1$ for some $P \in V(G)$.
\end{lemma}

\medskip

We now come to the main result of this section.  Returning to our
conventional notation, we let $R$ be a complete discrete valuation ring 
with field of fractions $K$ and algebraically closed
residue field $k$.
We let $X$ be a smooth curve over $K$,
and let $\X / R$ be a strongly semistable regular model for $X$ with special fiber $\X_k$.  

\begin{lemma}[Specialization Lemma]
\label{SpecializationLemma}
For all $D \in \Div(X)$, we have
\[
r_G(\rho(D)) \geq r_X(D).
\]
\end{lemma}

\begin{proof}
Let $\Dbar := \rho(D)$. 
We prove by induction on $k$ that if $r_X(D) \geq k$, then $r_G(\Dbar) \geq k$ as well.
The base case $k = -1$ is obvious.  Now suppose $k = 0$, so that $r_X(D) \geq 0$.
Then there exists an effective divisor $E \in \Div(X)$ with $D \sim E$, so
that $D - E \in \Prin(X)$.  Since $\rho$ is a homomorphism and 
takes principal (resp. effective) divisors on $X$ to principal (resp. effective) 
divisors on $G$, we have $\Dbar = \rho(D) \sim \rho(E) \geq 0$, so that 
$r_G(\Dbar) \geq 0$ as well.

We may therefore assume that $k \geq 1$.
Let $\Pbar \in V(G)$ be arbitrary.  By Remark~\ref{SurjectivityRemark},
there exists $P \in X(K)$ such that $\rho(P) = \Pbar$.
Then $r_X(D - P) \geq k-1$, so by induction we have $r_G(\Dbar - \Pbar) \geq k-1$ as well
(and in particular, $r_G(\Dbar) \geq 0$).
Since this is true for all $\Pbar \in V(G)$, it follows from
Lemma~\ref{GraphGenericPointLemma} that $r_G(\Dbar) \geq k$ as desired.
\end{proof}

\subsection{Compatibility with base change}

Let $K'/K$ be a finite extension, let $R'$ be the valuation ring of $K'$, 
and let $X_{K'} := X \times_K K'$.
It is known that there is a unique {\em relatively minimal} regular semistable 
model $\X'/R'$ which dominates $\X$, and the dual graph $G'$ of the special fiber
of $\X'$ is isomorphic to $\sigma_e(G)$, where $e$ is the ramification index of $K'/K$.
If we assign a length of $1/e$ to each edge of $G'$, then the corresponding {\em metric
graph} is the same for all finite extensions $K'/K$.  
In other words, $G$ and $G'$ are different models for the same metric $\QQ$-graph $\Gamma$, which we call
the {\em reduction graph} associated to the model $\X/R$ (see \cite{CR} for further discussion).

\medskip

The discussion in \cite{CR}, together with (\ref{ConcreteRho}), 
shows that
there is a unique surjective map 
$\tau : X(\Kbar) \to \Gamma_{\QQ}$ for which the induced homomorphism
\[
\tau_* : \Div(X_{\Kbar}) \cong \Div(X(\Kbar)) \to \Div_{\QQ}(\Gamma)
\]
is compatible with $\rho$, in the sense that 
for $D \in \Div(X(K'))$, we have $\tau_*(D) = \rho_{K'}(D)$.

Concretely, if $K'/K$ is a finite extension and $P \in X(K')$, then $\tau(P)$ is 
the point of $\Gamma_{\QQ}$ corresponding to the irreducible component of the special fiber of 
$\X'$ to which $P$ specializes.

\begin{remark}
In general, for $D \in \Div(X) \subseteq \Div(X_{\Kbar})$ we will not always have $\rho(D) = \tau_*(D)$,
but $\rho(D)$ and $\tau_*(D)$ will at least be linearly equivalent as divisors on the metric graph $\Gamma$.
(This is a consequence of standard facts from arithmetic intersection theory, e.g.,
Propositions 9.2.15 and 9.2.23 of \cite{Liu}.)
\end{remark}

From the discussion in \S\ref{MetricGraphSection}, we deduce from the proof of Lemma~\ref{SpecializationLemma}:

\begin{corollary}
\label{SpecializationBaseChangeCor}
Let $D \in \Div(X_{\Kbar})$.  Then
\[
r_{\QQ}(\tau_*(D)) \geq r_X(D).
\]
\end{corollary}

\begin{remark}
For the reader familiar with Berkovich's theory of analytic spaces \cite{BerkovichBook}, 
the following remark may be helpful.  
The metric $\QQ$-graph $\Gamma$ can be identified with the {\em skeleton} of the formal model
associated to $\X$, and the map $\tau : X(\Kbar) \to \Gamma_{\QQ}$ can be identified with the restriction to 
$X(\Kbar) \subset X^{\an}$ of the natural deformation retraction $X^{\an} \onto \Gamma$, where 
$X^{\an}$ denotes the Berkovich $K$-analytic space associated to $X$.
\end{remark}


\section{Some applications of the specialization lemma}
\label{ApplicationSection}

\subsection{Specialization of $g^r_d$'s}

Recall that a {\em complete $g^r_d$} on $X/K$ is defined to be a complete linear system
$|D|$ with $D \in \Div(X_{\Kbar})$, $\deg(D) = d$, and $r(D) = r$.
For simplicity, we will omit the word ``complete'' and just refer to such a linear system as
a $g^r_d$.
A $g^r_d$ is called {\em $K$-rational} if we can choose the divisor $D$ to lie in $\Div(X)$.

\medskip

By analogy, we define a $g^r_d$ on a graph $G$ (resp. a metric graph $\Gamma$) to be a 
complete linear system $|D|$ with $D \in \Div(G)$ (resp. $D \in \Div(\Gamma)$) such that 
$\deg(D) = d$ and $r(D) = r$.  Also, we will denote by $g^r_{\leq d}$ (resp. $g^{\geq r}_d$) a complete linear system $|D|$ 
with $\deg(D) \leq d$ and $r(D) = r$ (resp. $\deg(D) = d$ and $r(D) \geq r$).

As an immediate consequence of Lemma~\ref{SpecializationLemma} and Corollary~\ref{SpecializationBaseChangeCor}, 
we obtain:

\begin{corollary}
\label{grdcor}
Let $X$ be a smooth curve over $K$,
and let $\X / R$ be a strongly semistable regular model for $X$ with special fiber $\X_k$.  
Let $G$ be the dual graph of $\X_k$, and let $\Gamma$ be the corresponding metric graph.
If there exists a $K$-rational $g^r_d$ on $X$, then there exists a $g^{\geq r}_d$ and a $g^r_{\leq d}$ on $G$.
Similarly, if there exists a $g^r_d$ on $X$, then there exists a $g^{\geq r}_d$ 
and a $g^{r}_{\leq d}$ on $\Gamma$.
\end{corollary}

This result places restrictions on the possible graphs which can appear as the dual
graph of some regular model of a given curve $X/K$.

\medskip

In the particular case $r=1$, we refer to the smallest positive integer $d$ for which there exists a 
$g^1_d$ (resp. a $K$-rational $g^1_d$) on $X$ as the {\em gonality} (resp. {\em $K$-gonality}) of $X$.

Similarly, we define the {\em gonality} of a graph $G$ (resp. a metric graph $\Gamma$) to be the smallest
positive integer $d$ for which there exists a $g^1_d$ on $G$ (resp. $\Gamma$).

\medskip

As a special case of Corollary~\ref{grdcor}, we have:

\begin{corollary}
\label{gonalitycor}
The gonality of $G$ (resp. $\Gamma$) is at most the $K$-gonality (resp. gonality) of $X$.
\end{corollary}

\begin{example}
Let $K_n$ denote the complete graph on $n\geq 2$ vertices.

\medskip

{\bf Claim:} The gonality of $K_n$ is equal to $n-1$.

\medskip

Indeed, 
let $D = \sum a_v(v)$ be an effective divisor of degree at most $n-2$ on $K_n$, and
label the vertices $v_1,\ldots,v_n$ of $G$ so that $a_{v_1} \leq \cdots \leq a_{v_n}$.
Then it is easy to see that $a_{v_1} = 0$ and $a_{v_i} \leq i-2$ for all $2 \leq i \leq n$.
If $\nu$ is the divisor associated to the linear ordering 
$v_1 < \cdots < v_n$ of $V(K_n)$, it follows that $D - (v_1) \leq \nu$, 
so that $r(D - (v_1)) = -1$ by Theorem~\ref{OrderTheorem}.  In particular, we have 
$r(D) \leq 0$, and thus the gonality of $K_n$ is at least $n-1$.  
On the other hand, for any vertex $v_0 \in V(K_n)$, 
the divisor 
\[
D = \sum_{v \in V(K_n) \backslash \{ v_0 \} } (v)
\]
has degree $n-1$ and rank at least $1$, 
since $D - (v_0) \sim (n-2)(v_0)$.

\medskip

It follows from Corollary~\ref{gonalitycor} that if $X/K$ has $K$-gonality at most $n-2$, then
no regular model $\X/R$ for $X$ can have $K_n$ as its dual graph. 
For example, $K_4$ cannot be the dual graph of any regular model of a hyperelliptic curve $X/K$.
\end{example}

\subsection{Hyperelliptic graphs}

Focusing now on the special case $d=2$, we recall 
from \cite[\S{IV.5}]{Hartshorne} that a smooth curve $X/K$ of genus $g$ is 
called {\em hyperelliptic} if $g\geq 2$ and there exists a $g^1_2$ on $X$.  
If such a $g^1_2$ exists, it
is automatically unique and $K$-rational.  

Similarly, we say that a graph $G$ (resp. a metric graph $\Gamma$) of genus $g$ is {\em hyperelliptic} 
if $g \geq 2$ and there exists a $g^1_2$ on $G$ (resp. $\Gamma$).

\begin{remark}
One can show that if such a $g^1_2$ exists, it is automatically unique.
Also, if $G$ is 2-edge-connected of genus at least $2$, then $G$ is hyperelliptic if and only if 
there is an involution $h$ on $G$ for which the quotient graph $G/\la h \ra$ is a tree.  
These and other matters are discussed in \cite{BakerNorine2}.
\end{remark}

By Clifford's theorem for graphs \cite[Corollary 3.5]{BakerNorine}, 
if $g \geq 2$ and $D$ is a divisor of degree 2 on $G$ with 
$r(D) \geq 1$, then in fact $r(D) = 1$, and thus $G$ is hyperelliptic.
Combining this observation with Corollary~\ref{grdcor}, we find: 

\begin{corollary}
\label{HyperellipticCorollary}
If $X$ is hyperelliptic and $G$ has genus at least $2$, then $G$ is hyperelliptic as well.  
\end{corollary}

The converse of Corollary~\ref{HyperellipticCorollary} is false, as the following example shows.

\begin{example}
\label{BananaExample}

1. Let $G = B_n$ be the ``banana graph'' of genus $n-1$ consisting of 2 vertices $Q_1,Q_2$ connected by $n\geq 3$ edges.
Then the divisor $D = (Q_1)+(Q_2)$ on $G$ has degree 2 and $r(D) = 1$, so $G$ is hyperelliptic.
On the other hand, there are certainly non-hyperelliptic curves $X$ possessing a regular strongly semistable model 
with dual graph $G$.
For example, let $p \equiv 1 \pmod{12}$ be prime, and let $K = \QQ_p^{\nr}$ be the completion of the 
maximal unramified extension of $\QQ_p$.
Then the modular curve $X_0(p)$ has a regular semistable model over $K$ (the ``Deligne-Rapoport model''
\cite{DeligneRapoport}) whose dual graph is isomorphic to $B_n$ with $n = (p-1)/12$.  
However, by a result of Ogg \cite{OggHyperelliptic}, $X_0(p)$ is never hyperelliptic when $p>71$.

2. More generally, let $G = B(\ell_1,\ldots,\ell_n)$ be the graph obtained by subdividing the $i^{\rm th}$ edge of $B_n$
into $\ell_i$ edges for $1 \leq i \leq n$ (so that $B(1,1,\ldots,1) = B_n$).
Then one easily checks that $|(Q_1)+(Q_2)|$ is still a $g^1_2$, so that $G$ is hyperelliptic.
The dual graph of $\X$ is always of this type when the special fiber of $\X_k$ 
consists of two projective lines over $k$ intersecting transversely at $n$ points.
For example, the modular curve $X_0(p)$ with $p\geq 23$ prime has a regular model whose dual graph $G$ is of 
this type.
For all primes $p>71$, $G$ is hyperelliptic even though $X_0(p)$ is not.
\end{example}

\begin{remark}
Every graph of genus 2 is hyperelliptic, since by the Riemann-Roch theorem for graphs, the canonical
divisor $K_G$ has degree 2 and dimension 1.  
It is not hard to prove that for every integer $g \geq 3$, there are both hyperelliptic and non-hyperelliptic graphs 
of genus $g$.
\end{remark}

\subsection{Brill-Noether theory for graphs}

Classically, it is known by {\em Brill-Noether theory} 
that every smooth curve of genus $g$ over the complex numbers
has gonality at most $\lfloor (g+3)/2 \rfloor$, and this bound is tight: for every $g\geq 0$, 
the general curve of genus $g$ has gonality exactly equal to $\lfloor (g+3)/2 \rfloor$.
More generally:

\begin{theorem}[Classical Brill-Noether theory]
\label{BrillNoetherTheorem}
Fix integers $g,r,d \geq 0$, and define the {\em Brill-Noether number}
$\rho(g,r,d) = g - (r+1)(g-d+r)$.
Then:
\begin{enumerate}
\item If $\rho(g,r,d) \geq 0$, then every smooth curve $X/\CC$ 
of genus $g$ has a divisor $D$ with $r(D) = r$ and $\deg(D) \leq d$.
\item If $\rho(g,r,d)<0$, then on a general smooth curve of genus $g$, there is no 
divisor $D$ with $r(D) = r$ and $\deg(D) \leq d$.
\end{enumerate}
\end{theorem}

\medskip

Based on extensive computer calculations by Adam Tart (an undergraduate at Georgia Tech), 
we conjecture that similar results hold in the purely combinatorial setting of finite graphs:

\begin{conjecture}[Brill-Noether Conjecture for Graphs]
\label{BrillNoetherConj}
Fix integers $g,r,d \geq 0$, and set $\rho(g,r,d) = g - (r+1)(g-d+r)$.
Then:
\begin{enumerate}
\item If $\rho(g,r,d) \geq 0$, then every graph  
of genus $g$ has a divisor $D$ with $r(D) = r$ and $\deg(D) \leq d$.
\item If $\rho(g,r,d)<0$, there exists a graph of genus $g$ for which there is no 
divisor $D$ with $r(D) = r$ and $\deg(D) \leq d$.
\end{enumerate}
\end{conjecture}

In the special case $r=1$, Conjecture~\ref{BrillNoetherConj} can be reformulated as follows:

\begin{conjecture}[Gonality Conjecture for Graphs]
\label{GonalityConjecture}

For each integer $g \geq 0$:

(1) The gonality of any graph of genus $g$ is at most $\lfloor (g+3)/2 \rfloor$.

(2) There exists a graph of genus $g$ with gonality exactly $\lfloor (g+3)/2 \rfloor$.
\end{conjecture}

Adam Tart has verified part (2) of Conjecture~\ref{GonalityConjecture} for
$g \leq 12$, and part (2) of Conjecture~\ref{BrillNoetherConj} for $2 \leq r \leq 4$
and $g \leq 10$.  He has also verified that part (1) of Conjecture~\ref{GonalityConjecture}
holds for approximately $1000$ randomly generated graphs of genus at most $10$, and has similarly
verified part (1) of Conjecture~\ref{BrillNoetherConj} for around $100$ graphs in the case
$2 \leq r \leq 4$.

\medskip

Although we do not know how to handle the general case,
it is easy to prove that part (1) of Conjecture~\ref{GonalityConjecture} holds for small values of $g$:

\begin{lemma}
Part (1) of Conjecture~\ref{GonalityConjecture} is true for $g \leq 3$.
\end{lemma}

\begin{proof}
For $g \leq 2$, this is a straightforward consequence of Riemann-Roch for graphs.  
For $g=3$, we argue as follows.  The canonical divisor $K_G$ on 
any genus 3 graph $G$ has degree $4$ and $r(K_G) = 2$.  
By Lemma~\ref{GraphGenericPointLemma}, there exists a vertex
$P \in V(G)$ for which the degree $3$ divisor $K_G - P$ has rank 1.  (In fact, it is not
hard to see that $r(K_G - P)=1$ for {\em every} vertex $P$.)  Therefore $G$ has a $g^1_3$, so
that the gonality of $G$ is at most 3, proving the lemma.
\end{proof}

For {\em metric $\QQ$-graphs},
we can prove the analogue of part (1) of Conjecture~\ref{BrillNoetherConj}
using the Specialization Lemma and classical Brill-Noether theory:

\begin{theorem}
\label{GraphBrillNoetherTheorem}
Fix nonnegative integers $g,r,d$ for which $g - (r+1)(g-d+r) \geq 0$.
Then every metric $\QQ$-graph $\Gamma$ of genus $g$ has a divisor $D$ with $r(D) = r$ and $\deg(D) \leq d$.
\end{theorem}

\begin{proof}
By uniformly rescaling the edges of a suitable model for $\Gamma$ so that they all have integer lengths, 
then adding vertices of valence 2 as necessary,
we may assume that $\Gamma$ is the metric graph associated to a graph $G$
(and that every edge of $G$ has length 1). 
By Theorem~\ref{DeformationTheorem}, there exists an arithmetic surface $\X/R$ whose generic fiber is smooth 
of genus $g$ and whose special fiber has dual graph $G$.
By classical Brill-Noether theory, there exists a $g^r_{d'}$ on $X$
for some $d' \leq d$, so according to 
(the proof of) Corollary~\ref{grdcor},
there is a $\QQ$-divisor $D$ on $\Gamma$ with $\deg(D) \leq d$ and $r(D) = r$.
\end{proof}

In \S\ref{TropicalBrillNoetherSection}, we will generalize Theorem~\ref{GraphBrillNoetherTheorem} to arbitrary
metric graphs, and then to tropical curves, using ideas from \cite{GathmannKerber}.

\begin{remark}
It would be very interesting
to give a direct combinatorial proof of Theorem~\ref{GraphBrillNoetherTheorem}.  
In any case, we view the above proof of Theorem~\ref{GraphBrillNoetherTheorem} as an 
example of how one can use the Specialization Lemma, in conjunction with known theorems about
algebraic curves, to prove nontrivial results about graphs (or more precisely, in this case, 
{\em metric graphs}).
\end{remark}

For a given graph $G$ (resp. metric graph $\Gamma$) and an integer $r \geq 1$, let $D(G,r)$ (resp. $D(\Gamma,r)$) 
be the minimal degree $d$ of a $g^r_d$ on $G$ (resp. $\Gamma$). 

\begin{conjecture}
\label{SubdivisionConjecture}
Let $G$ be a graph, and let $\Gamma$ be the associated $\QQ$-graph.
Then for every $r \geq 1$, we have:
\begin{enumerate}
\item $D(G,r) = D(\sigma_k(G),r)$ for all $k \geq 1$.
\item $D(G,r) = D(\Gamma,r)$.
\end{enumerate}
\end{conjecture}

Adam Tart has verified part (1) of Conjecture~\ref{SubdivisionConjecture} for $100$ different graphs 
with $1 \leq r,k \leq 4$, and for $1000$ randomly generated graphs of genus up to $10$ 
in the special case where $r = 1$ and $k = 2$ or $3$.

Note that Conjecture~\ref{SubdivisionConjecture}, in conjunction with Theorem~\ref{GraphBrillNoetherTheorem},
would imply part (1) of Conjecture~\ref{BrillNoetherConj}.

\medskip

Finally, we have the following analogue of part (2) of Conjecture~\ref{BrillNoetherConj} for 
metric graphs:

\begin{conjecture}[Brill-Noether Conjecture for Metric Graphs]
\label{MetricBrillNoetherConj}
Fix integers $g,r,d \geq 0$, and set $\rho(g,r,d) = g - (r+1)(g-d+r)$.
If $\rho(g,r,d)<0$, then there exists a metric graph of genus $g$ for which there is no 
divisor $D$ with $r(D) = r$ and $\deg(D) \leq d$.
\end{conjecture}

Note that 
Conjecture~\ref{MetricBrillNoetherConj} would follow from Conjecture~\ref{SubdivisionConjecture}
and part (2) of Conjecture~\ref{BrillNoetherConj}.

\begin{remark}
By a simple argument based on Theorem~\ref{DeformationTheorem} and Corollary~\ref{gonalitycor},
a direct combinatorial proof of Conjecture~\ref{MetricBrillNoetherConj} in the special case $r=1$ 
would yield a new proof of the classical fact that for every $g \geq 0$, 
there exists a smooth curve $X$ of genus $g$ over an algebraically
closed field of characteristic zero having gonality at least $\lfloor (g+3)/2 \rfloor$.  
\end{remark}

\subsection{A tropical Brill-Noether theorem}
\label{TropicalBrillNoetherSection}

In this section, we show how the ideas from \cite{GathmannKerber} can be used to generalize
Theorem~\ref{GraphBrillNoetherTheorem} from metric $\QQ$-graphs to arbitrary metric graphs and 
tropical curves.

The key result is the following ``Semicontinuity Lemma'', which allows one to transfer certain results about divisors on $\QQ$-graphs 
to arbitrary metric graphs.  For the statement, fix a metric graph $\Gamma$ and a positive real number $\epsilon$ smaller than all 
edge lengths in some fixed model for $\Gamma$.  We denote by $A_\epsilon(\Gamma)$ the ``moduli space'' of all metric graphs that are
of the same combinatorial type as $\Gamma$, and whose edge lengths are within $\epsilon$ of the corresponding edge lengths in $\Gamma$.
Then $A_\epsilon(\Gamma)$ can naturally be viewed as a compact polyhedral complex; in particular,
there is a well-defined notion of convergence in $A_\epsilon(\Gamma)$.

Similarly, for each positive integer $d$, we define $M = M^{d}_{\epsilon}(\Gamma)$ to be the compact polyhedral complex
whose underlying point set is
\[
M := \{ (\Gamma',D') \; : \; \Gamma' \in A_\epsilon(\Gamma), \, D' \in \Div_+^d (\Gamma') \}.
\]

\begin{lemma}[Semicontinuity Lemma]
\label{SemicontinuityLemma}
The function $r : M \to \ZZ$ given by $r(\Gamma',D') = r_{\Gamma'}(D')$ is upper semicontinuous, i.e.,
the set $\{ (\Gamma',D') \; : \; r_{\Gamma'}(D') \geq i \}$ is closed for all $i$.
\end{lemma}

\begin{proof}
Following the general strategy of \cite[Proof of Proposition~3.1]{GathmannKerber}, but with some slight variations in notation,
we set
\[
S := \{ (\Gamma',D',f,P_1,\ldots,P_d) \; : \; \Gamma' \in A_\epsilon(\Gamma), \, D' \in \Div_+(\Gamma'),
\, f \in \M(\Gamma'), \, 
(f) + D' = P_1 + \cdots + P_d \}. 
\]
Also, for each $i = 0,\ldots,d$, set
\[
M_i := \{ (\Gamma',D',P_1,\ldots,P_i) \; : \; \Gamma' \in A_\epsilon(\Gamma), \, D' \in \Div_+^d(\Gamma'), \, 
P_1,\ldots,P_i \in \Gamma' \}. 
\]
As in Lemma~1.9 and Proposition 3.1 of \cite{GathmannKerber}, one can endow each of the spaces $S$ and $M_i$ ($0 \leq i \leq d$) 
with the structure of a polyhedral complex.  

The obvious ``forgetful morphisms'' 
\[
\pi_i : S \to M_i, \, \, (\Gamma',D',f,P_1,\ldots,P_d) \mapsto (\Gamma',D',P_1,\ldots,P_i)
\]
and 
\[
p_i : M_i \to M, \, \, (\Gamma',D',P_1,\ldots,P_i) \mapsto (\Gamma',D')
\]
are morphisms of polyhedral complexes, and in particular they are continuous maps between topological spaces.
Following \cite[Proof of Proposition~3.1]{GathmannKerber}, we make the following observations:
\begin{enumerate}
\item $p_i$ is an {\em open map} for all $i$ (since it is locally just a linear projection).
\item $M_i \backslash \pi_i(S)$ is a union of open polyhedra, and in particular, is an open subset of $M_i$.
\item For $(\Gamma',D') \in M$, we have $r_{\Gamma'}(D') \geq i$ if and only if $(\Gamma',D') \not\in p_i(M_i \backslash \pi_i(S))$.
\end{enumerate}

From (1) and (2), it follows that $p_i(M_i \backslash \pi_i(S))$ is open in $M$.
So by (3), we see that the subset $\{ (\Gamma',D') \; : \; r_{\Gamma'}(D') \geq i \}$ is closed in $M$, as desired.
\end{proof}

The following corollary shows that the condition for a metric graph to have a $g^r_{\leq d}$ is closed:

\begin{corollary}
\label{SemicontinuityCor}
Suppose $\Gamma_n$ is a sequence of metric graphs in $A_{\epsilon}(\Gamma)$ converging to $\Gamma$.
If there exists a $g^r_{\leq d}$ on $\Gamma_n$ for all $n$, then there exists a $g^r_{\leq d}$ on $\Gamma$ as well.
\end{corollary}

\begin{proof}
Without loss of generality, we may assume that $r \geq 0$.
Passing to a subsequence and replacing $d$ by some $d'\leq d$ if neccesary, we may assume that for each $n$, there exists
an effective divisor $D_n \in \Div_+(\Gamma_n)$ with $\deg(D_n) = d$ and $r(D_n) = r$.
Since $M$ is compact, $\{ (\Gamma_n, D_n) \}$ has a convergent subsequence; by passing to this subsequence, we may assume that 
$(\Gamma_n,D_n) \to (\Gamma,D)$ for some divisor $D \in \Div(\Gamma)$.
By Lemma~\ref{SemicontinuityLemma}, we have $r(D) \geq r$.  
Subtracting points from $D$ if necessary, we find that there is an effective divisor $D' \in \Div(\Gamma)$ with $\deg(D') \leq d$
and $r(D) = r$, as desired.
\end{proof}

\begin{corollary}
\label{SemicontinuityCor2}
Fix integers $r,g,d$.  If there exists a $g^r_{\leq d}$ on every $\QQ$-graph of genus $g$, then there exists a $g^r_{\leq d}$ on 
every metric graph of genus $g$.
\end{corollary}

\begin{proof}
We can approximate a metric graph $\Gamma$ by a sequence of $\QQ$-graphs in $A_{\epsilon}(\Gamma)$
for some $\epsilon > 0$, so the result follows directly from Corollary~\ref{SemicontinuityCor}.
\end{proof}

Finally, we give our promised application of the Semicontinuity Lemma to Brill-Noether theory for tropical curves:

\begin{theorem}
\label{TropicalBrillNoetherTheorem}
Fix integers $g,r,d \geq 0$ such that $\rho(g,r,d) = g - (r+1)(g-d+r) \geq 0$.
Then every tropical curve of genus $g$ has a divisor $D$ with $r(D) = r$ and $\deg(D) \leq d$.
\end{theorem}

\begin{proof}
By Theorem~\ref{GraphBrillNoetherTheorem}, there exists a $g^r_{\leq d}$ on every metric $\QQ$-graph, so
it follows from Corollary~\ref{SemicontinuityCor2} that the same is true for {\em all} metric graphs.
By \cite[Remark~3.6]{GathmannKerber}, if $\tilde{\Gamma}$ is a tropical curve and $\Gamma$ is the metric
graph obtained from $\tilde{\Gamma}$ by removing all unbounded edges, then for every $D \in \Div(\Gamma)$ we have
$r_{\tilde{\Gamma}}(D) = r_{\Gamma}(D)$.  Therefore the existence of a $g^r_{\leq d}$ on $\Gamma$ implies the 
existence of a $g^r_{\leq d}$ on $\tilde{\Gamma}$.
\end{proof}

\section{Weierstrass points on curves and graphs}
\label{WeierstrassSection}

As another illustration of the Specialization Lemma in action, in this section we will explore the relationship
between Weierstrass points on curves and (a suitable notion of) Weierstrass points on graphs.
As an application, we will generalize and place into a more conceptual framework a well-known result of Ogg 
concerning Weierstrass points on the modular curve $X_0(p)$.  We will also prove the existence of Weierstrass points
on tropical curves of genus $g \geq 2$.

\subsection{Weierstrass points on graphs}
\label{GraphWPSection}

Let $G$ be a graph of genus $g$.
By analogy with the theory of algebraic curves, we say that $P \in V(G)$ is a {\em Weierstrass point} if 
$r(g(P)) \geq 1$.
We define Weierstrass points on {\em metric graphs} and {\em tropical curves} in exactly the same way.

\begin{remark}
\label{WeierstrassSubdivisionRemark}
By Corollary~\ref{SubdivisionCorollary}, if $\Gamma$ is the $\QQ$-graph corresponding to $G$ (so that every edge
of $G$ has length 1), then $P \in V(G)$ is a Weierstrass point on $G$ if and only if $P$ is a Weierstrass point
on $\Gamma$.
\end{remark}

Let $P \in V(G)$.  An integer $k \geq 1$ is called a {\em Weierstrass gap} for $P$ if $r(k(P)) = r((k-1)(P))$.
The Riemann-Roch theorem for graphs, together with the usual arguments from the theory of algebraic curves, yields 
the following result, whose proof we leave to the reader:

\begin{lemma}
\begin{enumerate}
\item The following are equivalent:

\begin{enumerate}
\item $P$ is a Weierstrass point.
\item There exists a positive integer $k\leq g$ which is a Weierstrass gap for $P$.
\item $r(K_G - g(P)) \geq 0$.
\end{enumerate}

\item The set $W(P)$ of Weierstrass gaps for $P$ is contained in the finite set 
$\{ 1,2,\ldots,2g-1 \} \subseteq \NN$ and has cardinality $g$.
\end{enumerate}
\end{lemma}

\begin{remark}
\label{WeierstrassBananaRemark}

1. Unlike the situation for algebraic curves, there exist graphs of genus at least $2$ with no Weierstrass points.
For example, consider the graph $G = B_n$ of genus $g = n-1$ introduced in Example~\ref{BananaExample}.  We claim that
$B_n$ has no Weierstrass points for any $n\geq 3$.  Indeed, the canonical divisor $K_G$ is $(g-1)(Q_1) + (g-1)(Q_2)$, and by
symmetry it suffices to show that $r((g-1)(Q_2) - (Q_1)) = -1$.  This follows directly from Theorem~\ref{OrderTheorem},
since $(g-1)(Q_2) - (Q_1) \leq \nu := g(Q_2) - (Q_1)$ and $\nu$ 
is the divisor associated to the linear ordering $Q_1 < Q_2$ of $V(G)$.

2. More generally, let $G = B(\ell_1,\ldots,\ell_n)$ be the graph of genus $g = n-1$ 
obtained by subdividing the $i^{\rm th}$ edge of $B_n$ into $\ell_i$ edges.  
Let $R_{ij}$ for $1 \leq j \leq \ell_i$ denote the vertices strictly between $Q_1$ and $Q_2$ lying on the $i^{\rm th}$ edge
(in sequential order).  Then $Q_1$ and $Q_2$ are not Weierstrass points of $G$.  Indeed, by symmetry it again suffices to
show that $r((g-1)(Q_2) - (Q_1)) = -1$, and this follows from Theorem~\ref{OrderTheorem} by considering the linear
ordering 
\[
Q_1 < R_{11} < R_{12} < \ldots < R_{1(\ell_1 - 1)} < R_{21} < \cdots < R_{\ell_n (\ell_{n} - 1)} < Q_2.
\]

3. Other examples of families of graphs with no Weierstrass points are given in \cite{BakerNorine2}.
\end{remark}

\begin{example}
On the complete graph $G = K_n$ on $n\geq 4$ vertices, every vertex is a Weierstrass point.
Indeed, if $P,Q \in V(G)$ are arbitrary, then $g(P) - (Q)$ is equivalent to the effective divisor 
$(g - (n-1))(P) - (Q) + \sum_{v \in V(G)} (v)$, and thus $r(g(P)) \geq 1$. 
\end{example}

The following example, due to Serguei Norine, shows that there exist metric $\QQ$-graphs with infinitely many
Weierstrass points:

\begin{example}
\label{NorineExample}

Let $\Gamma$ be the metric $\QQ$-graph associated to the banana graph $B_n$ for some $n \geq 4$.
Then $\Gamma$ has infinitely many Weierstrass points.  

Indeed, label the edges of $\Gamma$ as
$e_1,\ldots,e_n$, and identify each $e_i$ with the segment $[0,1]$, where $Q_1$ corresponds to $0$, say, and
$Q_2$ corresponds to $1$.  We write $x(P)$ for the element of $[0,1]$ corresponding to the point $P \in e_i$
under this parametrization.  Then for each $i$ and each $P \in e_i$ with $x(P) \in [\frac{1}{3},\frac{2}{3}]$, 
we claim that $r(3(P)) \geq 1$, and hence $P$ is a Weierstrass point on $\Gamma$.

To see this, we will show explicitly that for every $Q \in \Gamma$ we have $|3(P) - (Q)| \neq \emptyset$.
For this, it suffices to construct a function $f \in \M(\Gamma)$ for which $\Delta(f) \geq -3(P) + (Q)$.
This is easy if $P = Q$, and otherwise we have:

\medskip

{\bf Case 1(a):} If $Q \in e_i$ and $x(P) < x(Q)$, let $y = \frac{x(P)+x(Q)}{2}$ and 
take $f$ to be constant on $e_j$ for $j \neq i$, and on $e_i$ to have 
slope $-2$ on $[y,x(P)]$, slope $1$ on $[x(P),x(Q)]$, and slope $0$ elsewhere.

\medskip

{\bf Case 1(b):} If $Q \in e_i$ and $x(Q) < x(P)$, let $y = \frac{3x(P)-x(Q)}{2}$ and 
take $f$ to be constant on $e_j$ for $j \neq i$, and on $e_i$ to have 
slope $-1$ on $[x(Q),x(P)]$, slope $2$ on $[x(P),y]$, and slope $0$ elsewhere.

\medskip

If $Q \in e_j$ for some $j \neq i$,
take $f$ to be constant on $e_k$ for $k \neq i,j$.
On $e_j$, let $z = \min(x(Q),1-x(Q))$, and 
take $f$ to have slope $1$ on $[0,z]$, slope $0$ on $[z,1-z]$, and slope $-1$ on
$[1-z,1]$.  Finally, along $e_i$, we have two cases:

\medskip

{\bf Case 2(a):} 
If $x(P) \in [\frac{1}{3},\frac{1}{2}]$,
let $y=3x(P)-1$, and define $f$ on $e_i$ to have 
slope $-1$ on $[0,y]$, slope $-2$ on $[y,x(P)]$, and slope $1$ on 
$[x(P),1]$.

\medskip

{\bf Case 2(b):} 
If $x(P) \in [\frac{1}{2},\frac{2}{3}]$,
let $y=3x(P)-1$, and define $f$ on $e_i$ to have 
slope $-1$ on $[0,x(P)]$, slope $2$ on $[x(P),y]$, and slope $1$ on 
$[y,1]$.
\end{example}

\begin{remark}
Similarly, one can show that for each integer $m \geq 2$, if $x(P) \in [\frac{1}{m},\frac{m-1}{m}]$ then
$r(m(P)) \geq 1$, and thus $P$ is a Weierstrass point on $\Gamma$ as long as the genus of $\Gamma$ is at least $m$.
\end{remark}

\medskip

We close this section with a result which generalizes Remark~\ref{WeierstrassBananaRemark}, and
which can be used in practice to identify non-Weierstrass points on certain graphs.

\begin{lemma}
\label{ResidualTreeLemma}
Let $v$ be a vertex of a graph $G$ of genus $g \geq 2$, and let $G'$ be the graph obtained by
deleting the vertex $v$ and all edges incident to $v$.  If $G'$ is a tree, 
then $v$ is not a Weierstrass point.
\end{lemma}

\begin{proof}
Since $G'$ is a tree, there is a linear ordering 
\[
v_1 < \cdots < v_{n-1}
\]
of the vertices of $G'$ such that each vertex
other than the first one has exactly one neighbor preceding it in the order.
Extend $<$ to a linear ordering of $V(G)$ by letting $v$ be the last element in the order.
Since the corresponding divisor $\nu$ is equal to $g(v) - (v_1)$, it follows from Theorem~\ref{OrderTheorem} that
$|g(v) - (v_1)| = \emptyset$, and therefore $r(g(v)) = 0$.  Thus $v$ is not a Weierstrass point.
\end{proof}

\begin{remark}
\label{ResidualTreeRemark}
It is easy to see that if $(G,v)$ satisfies the hypothesis of Lemma~\ref{ResidualTreeLemma}, then
so does $(\tilde{G},v)$, where $\tilde{G}$ is obtained by subdividing each edge $e_i$ of $G$ into $m_i$
edges for some positive integer $m_i$.
\end{remark}

\subsection{Specialization of Weierstrass points on totally degenerate curves}
\label{WeierstrassSpecializationSection}

We say that an arithmetic surface $\X/R$ is {\em totally degenerate} 
if the genus of its dual graph $G$ is the same as the genus of $X$.   
Under our hypotheses on $\X$, the genus of $X$ is the sum of the genus of $G$ and the genera of all irreducible components of $\X_k$,
so $\X$ is totally degenerate if and only if all irreducible components of $\X_k$ have genus 0.

\medskip

Applying the Specialization Lemma and the definition of a Weierstrass point, we immediately obtain:

\begin{corollary}
\label{WeierstrassSpecializationCor}
If $\X$ is a strongly semistable, regular, and totally degenerate arithmetic surface, 
then for every $K$-rational Weierstrass point $P \in X(K)$, 
$\rho(P)$ is a Weierstrass point of the dual graph $G$ of $\X$.  More generally, for every Weierstrass point
$P \in X(\Kbar)$, $\tau_*(P)$ is a Weierstrass point of the reduction graph $\Gamma$ of $\X$.
\end{corollary}

\medskip

As a sample consequence, we have the following concrete result:

\begin{corollary}
\label{WeierstrassBananaCor}
(1) Let $\X/R$ be a strongly semistable, regular, totally degenerate arithmetic surface whose
special fiber has a dual graph with no Weierstrass points (e.g. the graph 
$B_n$ for some $n \geq 3$).
Then $X$ does not possess any $K$-rational Weierstrass points.

(2) Let $\X/R$ be a (not necessarily regular) arithmetic surface whose
special fiber consists of two genus $0$ curves intersecting transversely at $3$ or more points.
Then every Weierstrass point of $X(K)$ specializes to a singular point of $\X_k$.

(3) More generally, let $\X/R$ be a (not necessarily regular) strongly semistable and totally degenerate 
arithmetic surface whose dual graph $G$ contains a vertex $v$ for which $G' := G \backslash \{ v \}$ 
is a tree.  Then there are no $K$-rational Weierstrass points on $X$ specializing to the component 
$C$ of $\X_k$ corresponding to $v$.
\end{corollary}

\begin{proof}
Part (1) follows from what we have already said.
For (2), it suffices to note that 
$X$ has a strongly semistable regular model $\X'$ whose 
dual graph $G'$ is isomorphic to $B(\ell_1,\ldots,\ell_g)$
for some positive integers $\ell_i$, and a point of $X(K)$ which specializes to a nonsingular point of $\X_k$ will specialize to 
either $Q_1$ or $Q_2$ in $G'$, neither of which is a Weierstrass point
by Remark~\ref{WeierstrassBananaRemark}.
Finally, (3) follows easily from Lemma~\ref{ResidualTreeLemma} and Remark~\ref{ResidualTreeRemark}.
\end{proof}

\medskip

We view Corollary~\ref{WeierstrassBananaCor} as a generalization of Ogg's argument in \cite{OggWeierstrass} showing 
that the cusp $\infty$ is never a Weierstrass point on $X_0(p)$ for $p \geq 23$ prime, since $X_0(p)/\QQ_p^{\nr}$ has a model $\X_0(p)$ of the type described in
part (2) of the corollary (the Deligne-Rapoport model), and $\infty$ specializes to a nonsingular point on the special fiber of $\X_0(p)$.
More generally, Corollary~\ref{WeierstrassBananaCor} shows (as does Ogg's original argument) that all Weierstrass points of
$X_0(p)$ specialize to supersingular points in the Deligne-Rapoport model.
Corollaries~\ref{WeierstrassSpecializationCor} and \ref{WeierstrassBananaCor}
give a recipe for extending Ogg's result to a much broader class of curves with 
totally degenerate reduction.

\begin{remark}
Corollary~\ref{WeierstrassBananaCor} has strong implications concerning the arithmetic of Weierstrass points on curves.
For example, in the special case where $K = \QQ_p^{\nr}$, 
the conclusion of part (1) of Corollary~\ref{WeierstrassBananaCor} implies that
every Weierstrass point on $X$ is ramified at $p$.
\end{remark}

\medskip

\begin{example}
The hypothesis that $\X$ is totally degenerate is necessary in the statement of 
Corollary~\ref{WeierstrassSpecializationCor}.
For example, it follows from a theorem of Atkin \cite[Theorem 1]{Atkin} 
that the cusp $\infty$ is a $\QQ$-rational Weierstrass point on $X_0(180)$.
The mod 5 reduction of the Deligne-Rapoport model of $X_0(180)$ consists of two copies of $X_0(36)_{\FF_5}$ intersecting transversely 
at the supersingular points, and the cusp $\infty$ specializes to a nonsingular point on one of these components.
(This does not contradict Corollary~\ref{WeierstrassBananaCor} because $X_0(36)$ does not have genus 0.)
\end{example}

\medskip

We conclude this section with another application of algebraic geometry to tropical geometry: we use
the classical fact that Weierstrass points exist on every smooth curve of genus at least 2 to show that 
there exist Weierstrass points on every tropical curve of genus at least 2.

\begin{theorem}
\label{WeierstrassExistenceTheorem}
Let $\tilde{\Gamma}$ be a tropical curve of genus $g \geq 2$.  
Then there exists at least one Weierstrass point on $\tilde{\Gamma}$.
\end{theorem}

\begin{proof}
We first consider the case of a $\QQ$-graph $\Gamma$.
By rescaling if necessary, we may assume that $\Gamma$ is the metric graph associated to a finite graph $G$,
with every edge of $G$ having length one.  
By Theorem~\ref{DeformationTheorem}, there exists a strongly semistable, regular, totally degenerate 
arithmetic surface $\X/R$ whose generic fiber is smooth 
of genus $g$ and whose special fiber has reduction graph $\Gamma$.
Let $P \in X(\Kbar)$ be a Weierstrass point, which exists by classical algebraic geometry since $g \geq 2$.  
By Corollary~\ref{WeierstrassSpecializationCor}, $\tau_*(P)$ is a Weierstrass point of $\Gamma$.
We have thus shown that every metric $\QQ$-graph of genus at least 2 has a Weierstrass point.

\medskip

Now let $\Gamma$ be an arbitrary metric graph.  As in \S\ref{TropicalBrillNoetherSection}, 
for $\epsilon > 0$ sufficiently small, we can approximate $\Gamma$ by a sequence $\Gamma_n$ of $\QQ$-graphs
within the space $A_{\epsilon}(\Gamma)$.  Let $P_n \in \Gamma_n$ be a Weierstrass point.  Passing to a subsequence
if necessary, we may assume without loss of generality that $(\Gamma_n,P_n) \to (\Gamma,P)$ in $M^{1}_\epsilon(\Gamma)$
for some point $P \in \Gamma$.  Since $r_{\Gamma_n}(gP_n) \geq 1$ for all $n$ and $(\Gamma_n,gP_n) \to (\Gamma,gP)$ in 
$M^{g}_{\epsilon}(\Gamma)$, we conclude from the Semicontinuity Lemma (Lemma~\ref{SemicontinuityLemma})
that $r_{\Gamma}(gP) \geq 1$, i.e., that $P$ is a Weierstrass point on $\Gamma$.

\medskip

Finally, suppose that $\tilde{\Gamma}$ is a tropical curve, and let 
$\Gamma$ be the metric graph obtained from $\tilde{\Gamma}$ by removing all unbounded edges.
It follows from \cite[Remark~3.6]{GathmannKerber} that every Weierstrass point on $\Gamma$ is also
a Weierstrass point on $\tilde{\Gamma}$.
Therefore the existence of Weierstrass points on $\Gamma$ implies the 
existence of Weierstrass points on $\tilde{\Gamma}$.
\end{proof}

\begin{remark}
We do not know a direct combinatorial proof of Theorem~\ref{WeierstrassExistenceTheorem}, but
it would certainly be interesting to give such a proof.
Also, in light of Example~\ref{NorineExample}, it is not clear if there is an analogue for
metric graphs of the classical fact that the total weight of all Weierstrass points on a smooth curve of
genus $g\geq 2$ is $g^3 - g$.
\end{remark}

\subsection{Specialization of a canonical divisor}


Since $P \in X(\Kbar)$ is a Weierstrass point of $X$ if and only if $r(K_X - g(P)) \geq 0$, where 
$K_X$ denotes a canonical divisor on $X$, and since
$P \in V(G)$ is a Weierstrass point of $G$ if and only if $r(K_G - g(P)) \geq 0$, 
Corollary~\ref{WeierstrassSpecializationCor} suggests a relationship between the canonical divisor of $G$ 
and the specialization of $K_X$ when $\X$ is totally degenerate.  
We investigate this relationship in this section.

\medskip

Let $\X$ be a strongly semistable, regular, totally degenerate arithmetic surface.
Let $\omega_{\X/R}$ be the {\em canonical sheaf} for $\X/R$, and let $K_{\X}$ be a Cartier
divisor such that $\O_{\X}(K_{\X}) \cong \omega_{\X/R}$; we call any such $K_{\X}$ a {\em canonical 
divisor}.  

\begin{lemma}
\label{CanonicalSpecializationLemma}
We have $\rho(K_{\X}) = K_G$.
\end{lemma}

\begin{proof}
This is a consequence of the {\em adjunction formula} for arithmetic surfaces (see \cite[Theorem 9.1.37]{Liu}), which
tells us that
\begin{equation}
\label{eq:Adjunction1}
(K_{\X} \cdot C_i) = 2g(C_i) - 2 - (C_i\cdot C_i) = -2 - (C_i \cdot C_i)
\end{equation}
for all $i$.
Since $(C \cdot C_i) = 0$ for all $i$, we have
\begin{equation}
\label{eq:Adjunction2}
(C_i \cdot C_i) = -\sum_{j \neq i} (C_i \cdot C_j) = -\deg(v_i).
\end{equation}

Combining (\ref{eq:Adjunction1}) and (\ref{eq:Adjunction2}) gives
\[
(K_{\X} \cdot C_i) = \deg(v_i) - 2 
\]
for all $i$, as desired.
\end{proof}

\begin{remark}
Lemma~\ref{CanonicalSpecializationLemma} helps explain why there is in fact a canonical {\em divisor} on a graph $G$,
rather than just a canonical {\em divisor class}, and also explains the connection between the canonical divisor on a graph
and the canonical divisor class in algebraic geometry.  This connection is implicit in the earlier work of S.~Zhang \cite{ZhangAP}.
\end{remark}

Note that the restriction of $K_{\X}$ to $X$ is a 
canonical divisor $K_X \in \Div(X)$ of degree $2g-2$, 
but $K_X$ is not necessarily supported on the set $X(K)$ of $K$-rational points of $X$.

\medskip

Let $\Gamma$ denote the metric $\QQ$-graph associated to $G$.
As before, we let $\tau : X(\Kbar) \to \Gamma_{\QQ}$ denote the natural surjective
specialization map, and we let $\tau_* : \Div(X(\Kbar)) \to \Div(\Gamma_{\QQ})$ denote the induced 
homomorphism on divisors.

\begin{lemma}
\label{CanonicalSpecializationLemma2}
Let $K_X \in \Div(X)$ be a canonical divisor.
Then $\tau_*(K_X)$ is linearly equivalent to $K_G$.
\end{lemma}

\begin{proof}
Since the Zariski closure of $K_X$ 
differs from a canonical divisor $K_{\X}$ on $\X$ by a vertical divisor, this follows from
Lemma~\ref{CanonicalSpecializationLemma}
and the remarks preceding Lemma~\ref{PrincipalSpecializationLemma}.
\end{proof}

\begin{remark}
\label{HorizontalRemark}
By a general {\em moving lemma} (e.g. Corollary 9.1.10 or Proposition 9.1.11
of \cite{Liu}), there exists a horizontal canonical divisor $K_{\X}$ on $\X$.  
Since $K_{\X}$ is the Zariski closure of $K_X$ in this case,
it follows that there exists a canonical divisor $K_X \in \Div(X(\Kbar))$ for which
$\tau_*(K_X)$ is {\em equal} to $K_G$ (and not just linearly equivalent to it).
\end{remark}

\subsection{An example}
\label{ExampleSection}

We conclude with an explicit example
which illustrates many of the concepts that have
been discussed in this paper.

\medskip

Let $p$ be an odd prime, and let $X$ be the smooth plane quartic curve over $\QQ_p$ given by
$F(x,y,z) = 0$, where
\begin{equation}
\label{eq:F}
F(x,y,z) = (x^2 - 2y^2 + z^2)(x^2 - z^2) + py^3z.
\end{equation}

By classical algebraic geometry, we have $g(X) = 3$, and the $\overline{\QQ}_p$-gonality of $X$ is
also $3$.  We let $K = \QQ_p^{\nr}$, and consider $X$ as an algebraic curve over $K$.

\medskip

Let $\X'$ be the model for $X$ over the valuation ring $R$ of $K$ given by the equation (\ref{eq:F}).
According to \cite[Exercise 10.3.10]{Liu}, the special fiber of $\X'$ is semistable and consists 
of two projective lines $\ell_1$ and $\ell_2$ with equations $x=z$ and $x=-z$, respectively, which intersect transversely
at the point $(0:1:0)$, together the conic $C$ defined by $x^2 - 2y^2 + z^2 = 0$, 
which intersects each of $\ell_1$ and $\ell_2$
transversely at 2 points.
The model $\X'$ is not regular, but a regular model $\X$ can be obtained from $\X'$ by blowing up the point $(0:1:0)$ of the
special fiber of $\X'$, which produces an exceptional divisor $E$ in $\X_k$ isomorphic to $\PP^1_k$, and which
intersects each of $\ell_1$ and $\ell_2$ transversely in a single point (see \cite[Corollary 10.3.25]{Liu}).
The special fiber $\X_k$ of $\X$ and the dual graph $G$ 
of $\X_k$ are depicted in Figures~\ref{reduction1} and \ref{reduction2} below.

\medskip

\begin{figure}[htb]
\begin{center}
\scalebox{.5}{\input{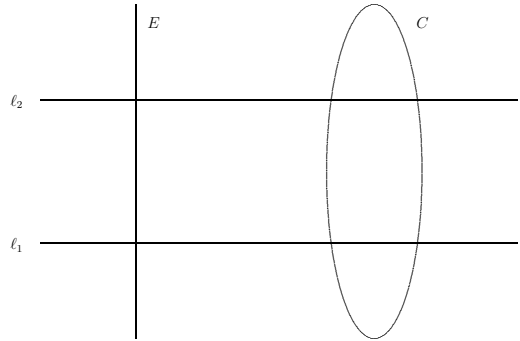}}  
\end{center}
\caption{The special fiber $\X_k$.}
\label{reduction1}
\end{figure}

\begin{figure}[htb]
\begin{center}
\scalebox{.5}{\input{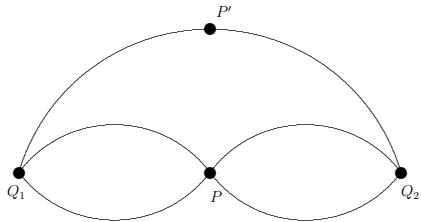}}  
\end{center}
\caption{The dual graph $G$ of $\X_k$.}
\label{reduction2}
\end{figure}

\medskip

In Figure~\ref{reduction2}, the vertex $P$ corresponds to the conic $C$, $Q_i$ corresponds to the line $\ell_i$ ($i=1,2$), 
and $P'$ corresponds to the exceptional divisor $E$ of the blowup.  Note that $G$ is a graph of
genus 3, and that $\X$ has totally degenerate strongly semistable reduction.

\medskip

{\bf Claim:} $Q_1$ and $Q_2$ are Weierstrass points of $G$, while $P$ and $P'$ are not.

\medskip

Indeed, since $3(Q_1) \sim 3(Q_2) \sim 2(P) + (P')$, it follows easily that 
$r(3(Q_1)) = r(3(Q_2)) \geq 1$, and therefore $Q_1$ and $Q_2$
are Weierstrass points.  On the other hand, $P$ is not a Weierstrass point of $G$ by 
Lemma~\ref{ResidualTreeLemma}, and $P'$ is not a Weierstrass point either, since 
$3(P') - (P)$ is equivalent to $(Q_1) + (Q_2) + (P') - (P)$, which is the divisor associated to
the linear ordering $P < Q_1 < Q_2 < P'$ of $V(G)$.

\medskip

{\bf Claim:} The gonality of $G$ is $3$.

\medskip

We have already seen that $r(3(Q_1)) \geq 1$, so the gonality of $G$ is at most 3.
It remains to show that $G$ is not hyperelliptic, i.e., that there is no $g^1_2$ on $G$.
By symmetry, and using the fact that $(Q_1) + (Q_2) \sim 2(P')$, it suffices to show that
$r(D) = 0$ for $D = (P)+(P')$ and for each of the divisors $(Q_1) + (X)$ with
$X \in V(G)$.
For this, it suffices to show that $|D| = \emptyset$ for 
each of the divisors $(Q_1) + 2(Q_2) - (P), (P) + 2(Q_1) - (P')$, and $(P) + (P') + (Q_1) - (Q_2)$.
But these are the $\nu$-divisors associated to the linear orderings 
$P < Q_1 < P' < Q_2$, $P' < Q_2 < P < Q_1$, and $Q_1 < P < Q_2 < P'$ of $V(G)$, respectively.
The claim therefore follows from Theorem~\ref{OrderTheorem}.

\medskip

The canonical divisor $K_G$ on $G$ is $2(P) + (Q_1) + (Q_2)$.
We now compute the specializations of various canonical divisors in 
$\Div(X(K))$.  Since a canonical divisor on $X$ is just a hyperplane section, the following divisors
are all canonical:

\begin{itemize}
\item[$K_1$:] $(x = z) \cap X = (0:1:0) + 3(1:0:1)$.
\item[$K_2$:] $(x = -z) \cap X = (0:1:0) + 3(1:0:-1)$.
\item[$K_3$:] $(z = 0) \cap X = 2(0:1:0) + (\sqrt{2}:1:0) + (-\sqrt{2}:1:0)$.
\item[$K_4$:] $(y = 0) \cap X = (1:0:1) + (1:0:-1) + (1:0:\sqrt{-1}) + (1:0:-\sqrt{-1})$.
\end{itemize}

The specializations of these divisors under $\rho$ (or equivalently, under $\tau_*$) 
are:

\begin{itemize}
\item[$\rho(K_1)$:] $(P') + 3(Q_1)$.
\item[$\rho(K_2)$:] $(P') + 3(Q_2)$.
\item[$\rho(K_3)$:] $2(P') + 2(P)$.
\item[$\rho(K_4)$:] $(Q_1) + (Q_2) + 2(P)$.
\end{itemize}

It is straightforward to check that each of these divisors is linearly equivalent to 
$K_G = (Q_1) + (Q_2) + 2(P)$,
in agreement with Lemma~\ref{CanonicalSpecializationLemma2}.

\medskip

Finally, note that (as follows from the above calculations) $(1:0:1)$ and $(1:0:-1)$ are Weierstrass points
on $X$, and they specialize 
to $Q_1$ and $Q_2$, respectively.  As we have seen, these are both Weierstrass 
points of $G$, as predicted by Corollary~\ref{WeierstrassSpecializationCor}.


\appendix
\section{A reformulation of Raynaud's description of the N{\'e}ron model of a Jacobian}
\label{RaynaudAppendix}

In this appendix, we re-interpret 
in the language of divisors on graphs some results
of Raynaud concerning the relation between a proper regular model for a curve
and the N{\'e}ron model of its Jacobian.  The main result here is
that the diagrams (\ref{BigDiagram1}) and (\ref{BigDiagram2}) below are exact and commutative.
This may not be a new observation, but since we could not find a reference, 
we will attempt to explain how it follows in a straightforward way from Raynaud's work.

In order to keep this appendix self-contained, we have repeated certain definitions which appear
in the main body of the text.
Some references for the results described here are 
\cite[Appendix]{BertoliniDarmon}, \cite{BLR}, \cite{Edixhoven}, and \cite{RaynaudPicard}.

\subsection{Raynaud's description}

Let $R$ be a complete discrete valuation ring with field of fractions $K$ and algebraically closed
residue field $k$.
Let $X$ be a smooth, proper, geometrically connected curve over $K$,
and let $\X / R$ be a proper model for $X$ with reduced special fiber $\X_k$.  
For simplicity, we assume throughout that $\X$ is regular, that 
the irreducible components of $\X_k$ are all
smooth, and that all singularities of $\X_k$ are ordinary double points.
We let $\C = \{ C_1,\ldots,C_n \}$ be the set of irreducible components of $\X_k$.

\medskip

Let $J$ be the Jacobian of $X$ over $K$, let $\J$ be the N{\'e}ron model of $J/R$,
and let $\J^0$ be the connected component of the identity in $\J$.
We denote by $\Phi = \J_k / \J_k^0$ the group of connected components of 
the special fiber $\J_k$ of $\J$.

\medskip

Let $\Div(X)$ (resp. $\Div(\X)$) be the group of Cartier divisors on $X$ (resp. on $\X$); 
since $X$ is smooth and $\X$ is regular, Cartier divisors on $X$ (resp. $\X$) are the same
as Weil divisors.  

The Zariski closure in $\X$ of an effective divisor on $X$ is a Cartier divisor.  Extending by
linearity, we can associate to each $D \in \Div(X)$ a Cartier divisor $\D$ on $\X$, which we
refer to as the {\em Zariski closure} of $D$.

Let $\Div^0(X)$ denote the subgroup of Cartier divisors of degree zero on $X$.
In addition, let $\Div^{(0)}(\X)$ denote the subgroup of $\Div(\X)$ consisting of those Cartier divisors $\D$ 
for which the restriction of the associated line bundle $\O_{\X}(\D)$
to each irreducible component of $\X_k$ has degree zero, i.e., 
for which 
\[
\deg(\O_{\X}(\D)|_{C_i}) = 0 \textrm{ for all } C_i \in \C.
\]

Finally, let
\[
\Div^{(0)}(X) = \{ D \in \Div^0(X) \; : \; \D \in \Div^{(0)}(\X) \},
\]
where $\D$ is the Zariski closure of $D$.

\medskip

Let $\Prin(X)$ (resp. $\Prin(\X)$) denote the group of principal Cartier divisors on $X$ (resp. $\X$).
There is a well-known isomorphism
\begin{equation}
\label{RaynaudIsom1}
J(K) = \J(R) \cong \Div^0(X) / \Prin(X),
\end{equation}
and according to Raynaud, there is an isomorphism
\begin{equation}
\label{RaynaudIsom2}
\begin{aligned}
J^0(K) := \J^0(R) \cong \Div^{(0)}(X) / \Prin^{(0)}(X),
\end{aligned}
\end{equation}
where 
\[
\Prin^{(0)}(X) := \Div^{(0)}(X) \cap \Prin(X).
\]

The isomorphism in (\ref{RaynaudIsom2}) 
comes from the fact that $\J^0 = \Pic^0_{\X/R}$ represents the functor
``isomorphism classes of line bundles whose restriction to each element of $\C$ has degree zero''.
(Recall that there is a canonical isomorphism between isomorphism classes of line bundles on $\X$ 
and the Cartier class group of $\X$.)

\begin{remark}
\label{KeyRemark}
In particular, it follows from the above discussion that every element $P \in J^0(K)$ can be
represented as the class of $D$ for some $D \in \Div^{(0)}(X)$.
\end{remark}

There is a natural inclusion $\C \subset \Div(\X)$, and 
an intersection pairing 
\[
\begin{aligned}
\C \times \Div(\X) &\to \ZZ \\
(C_i, \D) &\mapsto (C_i \cdot \D),
\end{aligned}
\]
where $(C_i \cdot \D) = \deg(\O_{\X}(D)|_{C_i})$.

\medskip

The intersection pairing gives rise to a map
\[
\begin{aligned}
\alpha : \ZZ^\C &\to \ZZ^\C \\
f &\mapsto \left( C_i \mapsto \sum_{C_j \in \C} (C_i \cdot C_j) f(C_j) \right).
\end{aligned}
\]

\medskip

Since $k$ is algebraically closed and the canonical map $J(K) \to \J_k(k)$ is surjective by
Proposition 10.1.40(b) of \cite{Liu},
there is a canonical isomorphism $J(K)/J^0(K) \cong \Phi$.
According to Raynaud, the component group $\Phi$ is canonically isomorphic to the homology of the complex
\begin{equation}
\label{RaynaudExactSeq1}
\begin{CD}
\ZZ^{\C} @>>{\alpha}> \ZZ^{\C} @>>{\deg}> \ZZ,
\end{CD}
\end{equation}
where $\deg : f \mapsto \sum_{C_i} f(C_i)$.

\medskip

The isomorphism 
\[
\phi : J(K)/J^0(K) \cong \Ker(\deg)/\Image(\alpha)
\]
can be described in the following way.  Let $P \in J(K)$, and choose $D \in \Div^0(X)$
such that $P = [D]$.  Let $\D \in \Div(\X)$ be the Zariski closure of $D$.
Then 
\[
\phi(P) = [C_i \mapsto (C_i \cdot \D)].
\]

\medskip

When $D$ corresponds to a Weil divisor supported on $X(K)$, we have another description
of the map $\phi$.  Write $D = \sum_{P \in X(K)} n_P(P)$ with $\sum n_P = 0$.
Since $\X$ is regular, each point $P \in X(K)=\X(R)$ specializes to a well-defined element $c(P)$ of $\C$.
Identifying a formal sum $\sum_{C_i \in \C} a_i C_i$ with the function $C_i \mapsto a_i \in \ZZ^{\C}$,
we have
\begin{equation}
\label{SpecializationSpecialCase}
\phi([D]) = [\sum_P n_P c(P)].
\end{equation}

\medskip

The quantities appearing in (\ref{RaynaudExactSeq1}) can be interpreted 
in a more suggestive fashion using the language of graphs.
Let $G$ be the {\em dual graph} of $\X_k$, i.e.,
$G$ is the finite graph whose vertices $v_i$ correspond to the irreducible components $C_i$ of $\X_k$,
and whose edges correspond to intersections between these components (so that there is one edge
between $v_i$ and $v_j$ for each point of intersection between $C_i$ and $C_j$).
We let $\Div(G)$ denote the free abelian group on the set of vertices of $G$,
and define $\Div^0(G)$ to be the kernel of the natural map $\deg : \Div(G) \to \ZZ$ given
by $\deg(\sum a_i (v_i)) = \sum a_i$.
In particular, the set $V(G)$ of vertices of $G$ is in bijection with $\C$, and
the group $\Div(G)$ is isomorphic to $\ZZ^{\C}$, with
$\Div^0(G)$ corresponding to $\Ker(\deg)$.

\medskip

Let $\M(G) = \ZZ^{V(G)}$ be the set of $\ZZ$-linear functions on $V(G)$, and
define the Laplacian operator
$\Delta : \M(G) \to \Div^0(G)$ by
\[
\Delta(\varphi) = \sum_{v \in V(G)} \sum_{e = vw} \left(\varphi(v) - \varphi(w) \right)(v),
\]
where the inner sum is over all edges $e$ of $G$ having $v$ as an endpoint.
Finally, define
\[
\Prin(G) = \Delta(\M(G)) \subseteq \Div^0(G),
\]
and let $\Jac(G) = \Div^0(G)/\Prin(G)$ be the {\em Jacobian} of $G$.  
It is a consequence of Kirchhoff's Matrix-Tree Theorem that $\Jac(G)$ is a finite abelian group whose order is equal to
the number of spanning trees of $G$.

\medskip

Since the graph $G$ is connected, one knows that $\Ker(\Delta)$ consists precisely of the
constant functions, and it follows from (\ref{RaynaudExactSeq1}) that 
there is a canonical exact sequence
\[
\begin{CD}
0 @>>> \Prin(G) @>\gamma_1>> \Div^0(G) @>\gamma_2>> \Phi @>>> 0.
\end{CD}
\]

In other words, the component group $\Phi$ is canonically isomorphic to the Jacobian group of the graph $G$.

\medskip

We can summarize much of the preceding discussion by saying that the following diagram is commutative and exact:

\begin{equation}
\label{BigDiagram1}
\begin{CD}
@.         0               @.            0          @.         0         @. \\
@.        @VVV                         @VVV                   @VVV       @. \\      
0 @>>> \Prin^{(0)}(X) @>\alpha_1>> \Div^{(0)}(X) @>\alpha_2>> J^0(K) @>>> 0  \\
@.        @VVV                         @VVV                   @VVV       @. \\      
0 @>>> \Prin(X)      @>\beta_1>>  \Div^{0}(X)   @>\beta_2>>  J(K)   @>>> 0  \\
@.        @VVV                         @VVV                   @VVV       @. \\ 
0 @>>> \Prin(G)      @>\gamma_1>>  \Div^{0}(G)   @>\gamma_2>> \Jac(G) \cong \Phi   @>>> 0  \\
@.        @VVV                         @VVV                   @VVV       @. \\
@.         0               @.            0          @.         0         @. \\      
\end{CD}
\end{equation}

A few remarks are in order about the exactness of the rows and columns in (\ref{BigDiagram1}).  
It is well-known that the natural map from $X(K)=\X(R)$ to the smooth locus of $\X_k(k)$ is surjective (see e.g.
Proposition 10.1.40(b) of \cite{Liu}); by
(\ref{SpecializationSpecialCase}), this implies that
the natural maps $\Div(X) \to \Div(G)$ and $\Div^0(X) \to \Div^0(G)$ are surjective.
The surjectivity of the horizontal map $\alpha_2 : \Div^{(0)}(X) \to J^0(K)$ follows from Remark~\ref{KeyRemark}.
Using this, we see from the Snake Lemma that since the vertical map $\Div^0(X) \to \Div^0(G)$ is surjective, 
the vertical map $\Prin(X) \to \Prin(G)$ is also surjective.
All of the other claims about the commutativity and exactness of
(\ref{BigDiagram1}) follow in a straightforward way from the definitions.

\subsection{Passage to the limit}

If $K'/K$ is a finite extension of degree $m$ with ramification index $e \mid m$ and valuation ring $R'$, 
then by a sequence of blow-ups 
we can obtain a regular model $\X'/R'$ for $X$ whose corresponding dual graph $G'$ is the graph $\sigma_e(G)$
obtained by subdividing each edge of $G$ into $e$ edges. 
If we think of $G$ as an unweighted graph and of $\sigma_e(G)$ 
as a weighted graph in which every edge has length $1/e$, then
$G$ and $\sigma_e(G)$ are different models for the same metric $\QQ$-graph $\Gamma$, which one calls
the {\em reduction graph} of $\X/R$.
The discussion in \cite{CR} shows that the various maps $c_{K'} : X(K') \to G'$ 
are compatible, in the sense that they give rise to a specialization map
$\tau : X(\Kbar) \to \Gamma$
which takes $X(\Kbar)$ surjectively onto $\Gamma_{\QQ}$.

It is straightforward to check that the diagram (\ref{BigDiagram1}) behaves 
functorially with respect to finite extensions, and therefore that there
is a commutative and exact diagram

\begin{equation}
\label{BigDiagram2}
\begin{CD}
@.         0               @.            0          @.         0         @. \\
@.        @VVV                         @VVV                   @VVV       @. \\      
0 @>>> \Prin^{(0)}(X(\Kbar)) @>\alpha_1>> \Div^{(0)}(X(\Kbar)) @>\alpha_2>> J^0(\Kbar) @>>> 0  \\
@.        @VVV                         @VVV                   @VVV       @. \\      
0 @>>> \Prin(X(\Kbar))      @>\beta_1>>  \Div^{0}(X(\Kbar))   @>\beta_2>>  J(\Kbar)   @>>> 0  \\
@.        @VVV                         @VVV                   @VVV       @. \\ 
0 @>>> \Prin_{\QQ}(\Gamma)      @>\gamma_1>>  \Div^{0}_{\QQ}(\Gamma)   @>\gamma_2>> \Jac_{\QQ}(\Gamma)   @>>> 0  \\
@.        @VVV                         @VVV                   @VVV       @. \\
@.         0               @.            0          @.         0         @. \\      
\end{CD}
\end{equation}

\begin{remark}
Let $\tK$ be the completion of an algebraic closure $\Kbar$ of $K$, so that 
$\tK$ is a complete and algebraically closed field equipped 
with a valuation $v : \tK \to \QQ \cup \{ +\infty \}$, and $\Kbar$ is dense in $\tK$.
By continuity, one can extend $\tau$ to a map $\tau : X(\tK) \to \Gamma$ 
and replace $\Kbar$ by $\tK$ everywhere in the diagram (\ref{BigDiagram2}).
\end{remark}

A few explanations are in order concerning the definitions of the various groups and 
group homomorphisms which appear in (\ref{BigDiagram2}).
Since $\Kbar$ is algebraically closed, we may identify the group $\Div(X_{\Kbar})$ of Cartier (or Weil) 
divisors on $X_{\Kbar}$ with $\Div(X(\Kbar))$, the free abelian group on the set $X(\Kbar)$.
We define $\Prin(X(\Kbar))$ to be the subgroup of $\Div(X(\Kbar))$ consisting of principal divisors.
The group $\Div(X(\Kbar))$ (resp. $\Prin(X(\Kbar))$) can be identified with the direct limit of
$\Div(X_{K'})$ (resp. $\Prin(X_{K'})$) over all finite extensions $K'/K$.
Accordingly, we define the group $J^0(\Kbar)$ 
to be the direct limit of the groups $J^0(K')$ 
over all finite extensions $K'/K$,
and we define $\Div^{(0)}(X(\Kbar))$ and $\Prin^{(0)}(X(\Kbar))$ similarly.
Finally, we define $\Jac_{\QQ}(\Gamma)$ to be the quotient
$\Div^{0}_{\QQ}(\Gamma) / \Prin_{\QQ}(\Gamma)$.

The fact that $\Prin_{\QQ}(\Gamma)$, as defined in \S\ref{IntroSection}, coincides with the direct limit over all finite extensions $K'/K$ 
of the groups $\Prin(G')$ follows easily from Remark~\ref{LaplacianRemark}. 

With these definitions in place, it is straightforward to check using 
(\ref{BigDiagram1}) that the diagram (\ref{BigDiagram2}) is both commutative and exact.

\medskip

In particular, we note the following consequence of the exactness of 
(\ref{BigDiagram1}) and (\ref{BigDiagram2}):

\begin{corollary}
\label{SurjectiveCor}
The canonical maps $\Prin(X) \to \Prin(G)$ and
$\Prin(X(\Kbar)) \to \Prin_{\QQ}(\Gamma)$ are surjective.
\end{corollary}

\begin{remark}
It follows from Corollary~\ref{SurjectiveCor} that if $G$ is a graph, the group $\Prin(G)$ can
be characterized as the image of $\Prin(X)$ under the specialization map from $\Div(X)$ to $\Div(G)$
for any regular arithmetic surface $\X/R$ whose special fiber has dual graph isomorphic
to $G$.  (Such an $\X$ always exists by Corollary~\ref{DeformationCor} below.)
\end{remark}

\begin{remark}
Another consequence of (\ref{BigDiagram2}) is that there is a canonical isomorphism
\[
J(\Kbar) / J^0(\Kbar) \cong \Jac_{\QQ}(\Gamma),
\]
so that the group $\Jac_{\QQ}(\Gamma)$ plays the role of the component group of the N{\'e}ron model in this situation,
even though there is not a well-defined N{\'e}ron model for $J$ over $\Kbar$ or $\tK$, since the valuations on 
these fields are not discrete.
One can show using elementary methods that $\Jac_{\QQ}(\Gamma)$ is (non-canonically) isomorphic
to $(\QQ/\ZZ)^g$ (compare with the discussion in \cite[Expos{\'e} IX, \S{11.8}]{SGA7}).
\end{remark}

\section{A result from the deformation theory of stable marked curves (written by Brian Conrad)}
\label{DeformationTheorySection}


Recall from \S\ref{section:Notation} that by an {\em arithmetic surface}, we mean
a proper flat scheme of dimension 2 over a discrete valuation ring whose generic fiber is a smooth curve.
In this appendix, we describe how one can realize an arbitrary graph $G$ as the dual
graph of the special fiber of some regular arithmetic surface whose special fiber
$C$ is a {\em totally degenerate semistable curve} (or {\em Mumford curve}),
meaning that $C$ is semistable, every irreducible
component of $C$ is isomorphic to the projective line over the residue
field $k$, and all singularities of $C$ are $k$-rational.

\medskip

We begin with the following simple lemma, whose proof is left to the reader:

\begin{lemma}
\label{DeformationLemma}
Let $G$ be a connected graph, and let $k$ be an infinite field.  Then
there exists a totally degenerate semistable curve $C/k$ whose dual graph is
isomorphic to $G$.
\end{lemma}

The crux of the matter is the following theorem, whose proof is a standard
application of the deformation theory of stable marked curves.

\begin{theorem}
\label{DeformationTheorem}  
Let $C$ be a proper and geometrically connected
semistable curve over a field $k$, and let $R$ be a complete
discrete valuation ring with residue field $k$.
Assume that $k$ is infinite.
Then there exists an arithmetic surface
$\mathfrak{X}$ over $R$ with special fiber $C$ 
such that $\mathfrak{X}$ is a regular scheme.
\end{theorem}

Combining these two results, we obtain:

\begin{corollary}
\label{DeformationCor}
Let $R$ be a complete discrete valuation ring with field of fractions $K$ and
infinite residue field $k$.  For any connected graph $G$, there exists a 
regular arithmetic surface $\X / R$ whose generic fiber is a smooth, proper, and 
geometrically connected curve $X/K$, and whose special fiber is a totally degenerate 
semistable curve with dual graph isomorphic to $G$.
\end{corollary}

\begin{remark}
By Proposition 10.1.51 of \cite{Liu}, the genus of $X$
coincides with the genus $g = |E(G)| - |V(G)| + 1$
of the graph $G$.
\end{remark}

\begin{proof}[Proof of Theorem~\ref{DeformationTheorem}]
We work with a general residue field $k$ (and avoid assumptions on
the field of rationality of the singularities) until near the end of
the argument. 
Let $g = \dim {\rm{H}}^1(C,\mathscr{O}_C)$ denote the arithmetic genus of $C$.
The structure theorem for ordinary double points
\cite[III,~\S2]{FreitagKiehl} ensures that $C^{\rm{sing}}$ splits over
a separable extension of $k$, so we can choose 
a finite Galois extension $k'/k$ 
so that the locus $C_{k'}^{\rm{sing}}$ of non-smooth points in $C_{k'}$ consists entirely of 
$k'$-rational points and every 
irreducible component of $C_{k'}$ is geometrically
irreducible and has a $k'$-rational point.  (If $C$ is a 
Mumford curve then we can take $k' = k$.)
In particular, each smooth component in $C_{k'}$ of arithmetic genus 0
is isomorphic to $\mathbf{P}^1_{k'}$ and so admits 
at least three $k'$-rational points.
We can then construct
a ${\rm{Gal}}(k'/k)$-stable \'etale divisor $D' \subseteq C_{k'}^{\rm{sm}}$ whose support consists
entirely of $k'$-rational points in the smooth locus such that for each component $X'$ of
$C_{k'}$ isomorphic to $\mathbf{P}^1_{k'}$ we have
$\#(X' \cap C_{k'}^{\rm{sing}}) +   \#(X' \cap D') \ge 3$.  In particular, if we choose
an enumeration of $D'(k')$ then the pair $(C_{k'},D')$ is a
stable $n$-pointed genus-$g$ curve, where $n = \#D'(k')$
and $2g-2 + n > 0$.  Let $D \subseteq C^{\rm{sm}}$ be the
\'etale divisor that descends $D'$.    We let $R'$ be the local finite
\'etale $R$-algebra with residue field $k'/k$. 

The stack $\mathscr{M}_{g,n}$ classifying stable $n$-pointed genus-$g$ curves
for any $g, n \ge 0$ such that $2g-2+n > 0$ 
is a proper smooth Deligne--Mumford stack over $\Spec \ZZ$.\footnote{
This is a standard fact:  it is due to Deligne--Mumford \cite[5.2]{DeligneMumford} if $n = 0$ (so $g \ge 2$),
Deligne--Rapoport if $g = n = 1$ \cite[IV,~2.2]{DeligneRapoport}, and is trivial if $g = 0, n = 3$ (in which
case $(\mathbf{P}^1, \{0, 1, \infty\})$ is the only such object, so the stack is
$\Spec \ZZ$).   The general case follows from these cases by realizing
$\mathscr{M}_{g,n}$ as the universal curve over $\mathscr{M}_{g,n-1}$
(due to Knudsen's
contraction and stabilization operations \cite[2.7]{Knudsen2}).}
The existence of $\mathscr{M}_{g,n}$ as a smooth Deligne--Mumford stack ensures that
$(C_{k'},D')$ admits a universal formal deformation 
$(\widehat{\mathscr{C}'}, \widehat{\mathscr{D}'})$ over 
a complete local noetherian $R'$-algebra $A'$ with residue field $k'$,
and that $A'$ is a formal power series ring over $R'$.
Moreover, there is a relatively ample
line bundle canonically associated to any stable $n$-pointed genus-$g$ curve $(X, \{\sigma_1,
\dots, \sigma_n\}) \rightarrow S$ over a scheme $S$, 
namely the twist $\omega_{X/S}(\sum \sigma_i)$ of the relative dualizing sheaf of the curve 
by the \'etale divisor defined by the marked points.   Hence, 
the universal formal deformation uniquely algebraizes to a pair
$(\mathscr{C}', \mathscr{D}')$ over $\Spec A'$.

Since $(C_{k'},D') = k' \otimes_k (C,D)$, 
universality provides an action of the Galois group $\Gamma = {\rm{Gal}}(k'/k)$ on $A'$
and on $(\mathscr{C}', \mathscr{D}')$ covering the natural $\Gamma$-action on $(C_{k'},D')$
and on $R'$. 
The action by $\Gamma$ on $\mathscr{C}'$ is compatible with one
on the canonically associated
ample line bundle $\omega_{\mathscr{C}'/A'}(\mathscr{D}')$.  Since $R \rightarrow R'$ is finite \'etale with
Galois group $\Gamma$, the $\Gamma$-action on everything in sight 
(including the relatively ample line bundle) defines effective descent data:
$A = (A')^{\Gamma}$ is a complete local noetherian $R$-algebra
with residue field $k$ such that $R' \otimes_R A \rightarrow A'$
is an isomorphism, and $(\mathscr{C}', \mathscr{D}')$
canonically descends to a deformation $(\mathscr{C}, \mathscr{D})$ of $(C,D)$ over $A$. 
This can likewise be shown to be a universal deformation of $(C,D)$
in the category of complete local noetherian $R$-algebras
with residue field $k$, but we do not need this fact.

What matters for our purposes is structural information about $A$ and
the proper flat $A$-curve $\mathscr{C}$.   By the functorial characterization of
formal smoothness, $A$ is formally smooth over $R$ 
since the same holds for $A'$ over $R'$, so $A$ is a power series
ring over $R$ in finitely many variables.  
We claim that the Zariski-open locus of smooth curves in $\mathscr{M}_{g,n}$ is 
dense, which is to say that every geometric point has a smooth deformation.
Indeed, since Knudsen's contraction and stabilization operations
do nothing to smooth curves, this claim immediately reduces to the special
cases $n = 0$ with $g \ge 2$, $g = n = 1$, and $g = 0, n = 3$.
The final two cases are obvious and the first case was proved
by Deligne and Mumford via deformation theory \cite[1.9]{DeligneMumford}.  It follows
that the generic fiber of $\mathscr{C}'$ over $A'$ is a smooth curve, so
the same holds for $\mathscr{C}$ over $A$.  In other words,
the Zariski-closed locus in $\mathscr{C}$ where
$\Omega^1_{\mathscr{C}/A}$ is not invertible has
its closed image in $\Spec(A)$ given by a closed
subset $\Spec(A/I)$ for a unique nonzero radical ideal $I \subseteq A$.  
Since $A$ is a power series ring over $R$ and $I$ is a nonzero 
ideal, we can certainly find a local $R$-algebra map $\phi:A \rightarrow R$
in which $I$ has nonzero image.  The pullback $\mathfrak{X}_{\phi}$ of $\mathscr{C}$
along $\phi$ is a proper flat semistable curve over $R$ deforming $C$
such that the generic fiber is smooth (as otherwise
the map $\Spec(R) \rightarrow \Spec(A)$ would factor
through $\Spec(A/I)$, a contradiction).    

The only remaining problem is to show that if $\phi$ is chosen 
more carefully then $\mathfrak{X}_{\phi}$ is also regular. 
If $C$ is $k$-smooth then $\mathscr{C}$ is $A$-smooth, so
$\mathfrak{X}_{\phi}$ is $R$-smooth (and in particular regular). 
Thus, we may assume that $C$ is not $k$-smooth.
The main point is to use an understanding of the structure of
$\mathscr{C}$ near each singular point $c \in C - C^{\rm{sm}}$.
We noted at the outset that each finite extension
$k(c)/k$ is separable, and if $R_c$ is the corresponding
local finite \'etale extension of $R$ then the structure theory of ordinary double points
\cite[III,~\S2]{FreitagKiehl} provides an $R_c$-algebra isomorphism 
$\mathscr{O}_{\mathscr{C},c}^{\rm{h}} \simeq
A_c[u,v]^{\rm{h}}/(uv - a_c)$ for
some nonzero non-unit $a_c \in A_c$, where
$A_c$ is the local finite \'etale extension of $A$
with residue field $k(c)/k$; we have $a_c \ne 0$ since
$\mathscr{C}$ has smooth generic fiber.    
For $B = \ZZ[t,u,v]/(uv - t)$, the $B$-module 
$\wedge^2(\Omega^1_{B/\ZZ[t]})$ is $B/(u,v) = B/(t)$.
Thus, by a formal computation at each $c$
we see that the annihilator ideal of $\wedge^2(\Omega^1_{\mathscr{C}/A})$
on $\mathscr{C}$ cuts out an $A$-finite closed subscheme of the non-smooth
locus in $\mathscr{C}$ that is a pullback of a unique
closed subscheme 
$\Spec(A/J) \subseteq \Spec(A)$.    We made the initial choice
of $k'/k$ large enough so that it splits each $k(c)/k$.
Hence,  the method of proof of \cite[1.5]{DeligneMumford} and
the discussion following that result show that 
for each singularity $c'$ in $C_{k'}$, the corresponding element
$a_{c'} \in A'$ may be chosen so that the $a_{c'}$'s are part
of a system of variables for $A'$ as a formal power series ring over $R'$.  
The ideal $J A'$ is the intersection of the ideals $(a_{c'})$.
To summarize, $A$ is a formal power series
ring over $R$ and we can choose the variables
for $A'$ over $R'$ such that $JA'$ is  generated by a product
of such variables, one for each singularity on $C_{k'}$.  
In particular, the local interpretation of each $a_{c'}$ on $\mathscr{C}'$
shows that for a local $R'$-algebra map $\phi':A' \rightarrow R'$, the pullback
$\mathfrak{X}'_{\phi'}$ of $\mathscr{C}'$ along $\phi'$ is regular if and only if
$\phi'(a_{c'}) \in R'$ is a uniformizer for each $c'$.   This condition on $\phi'$
is equivalent to saying that $\phi'(J A') \subseteq R'$ is
a proper nonzero ideal with multiplicity equal to the
number $\nu$ of 
geometric singularities on $C$.

Since $JA' = J \otimes_A A'$ is a principal nonzero proper ideal, and hence it is 
invertible as an $A'$-module, it follows that $J$ is principal as
well, say $J = (\alpha)$ for some nonzero nonunit $\alpha \in A$. 
We seek an $R$-algebra map $\phi:A \rightarrow R$ such that
$\phi(I) \ne 0$ and $\phi(\alpha) \in R$ is nonzero with order $\nu$, for then
$\phi^{\ast}(\mathscr{C})$ will have smooth generic
fiber (by our earlier discussion) and will become regular over $R'$, and so
it will be regular since $R'$ is finite \'etale over $R$.  The information
we have about $\alpha \in A$ is that in 
the formal power series ring $A'$ over $R'$ we can choose the variables so that 
$\alpha$ is a product of $\nu$ of the variables. 
We are now reduced to the following problem in
commutative algebra.  Let $R \rightarrow R'$ be a local finite \'etale
extension of discrete valuation rings, $A = R[\![x_1,\dots,x_N]\!]$,
$I \subseteq A$ a nonzero ideal,
and $\alpha \in A$ an element such that in $A' = R' \otimes_R A$
we can write $\alpha = x'_1 \cdots x'_{\nu}$ for some $1 \le \nu \le N$
and choice of $R'$-algebra isomorphism $A' \simeq R'[\![x'_1,\dots,x'_N]\!]$.
Then we claim that there is an $R$-algebra map
$\phi:A \rightarrow R$ such that ${\rm{ord}}_R(\phi(\alpha)) = \nu$
and $\phi(I) \ne 0$.

If we can choose $k' = k$, such as in the case when $C$ is a 
Mumford curve, then such a $\phi$ obviously exists.  Hence, for
the intended application to Corollary~\ref{DeformationCor}, we are done. 
To prove the claim in general, consider the expansions
$$x'_j = a'_{0j} + \sum_{i=1}^N a'_{ij} x_i + \dots$$
where $a'_{0j} \in \mathfrak{m}_{R'}$ and $(a'_{ij})_{1 \le i, j \le N}$
is invertible over $R'$.   Let $\pi \in \mathfrak{m}_R$ be a uniformizer.
We seek $\phi$ of the form $\phi(x_i) = t_i \pi$ for
$t_i \in R$.    The requirement on the $t_i$'s is
that $(a'_{0j}/\pi) + \sum a'_{ij} t_i \in R'^{\times}$
for $1 \le i \le \nu$ and that $h(t_1 \pi, \dots, t_N \pi) \ne 0$
for some fixed nonzero power series $h \in I$.  The unit conditions
only depend on $t_i \bmod \mathfrak{m}_R$.  Thus, once
we find $t_i$ that satisfy these unit conditions,
the remaining non-vanishing condition on the nonzero
power series $h$ is trivial to satisfy by modifying the higher-order
parts of the $t_i$'s appropriately.  It remains to consider the unit conditions,
which is a consequence of the trivial lemma below.
\end{proof}

\begin{lemma} Let $k'/k$ be a finite extension of fields.  
For $1 \le \nu \le N$ let $\{H'_1, \dots, H'_{\nu}\}$ be a collection
of independent hyperplanes in ${k'}^N$.  If $k$ is infinite then for any 
$v'_1, \dots, v'_{\nu} \in {k'}^N$, the union of the affine-linear
hyperplanes $v'_i + H'_i$ in ${k'}^N$ cannot contain $k^N$. 
\end{lemma}

If this lemma is true for finite $k$ then Theorem~\ref{DeformationTheorem} is
true without restriction on $k$.    By a long induction argument
that we omit, this lemma for a fixed finite field $k$ can be reduced to
establishing an affirmative answer to the following 
question in the theory of subspace arrangements over 
$k$:  if $V$ is a vector space of dimension $N \ge 1$
over a finite field $k$, and if $V_1, \dots, V_N$ are nonzero
proper linear subspaces such that $V$ is a union of
the translates $v_i + V_i$ for $v_1, \dots, v_N \in V$,
is there some $1 \le d \le N$ for which
some $d$-fold intersection $V_{i_1} \cap \dots \cap V_{i_d}$ has
dimension at least $N-d+1$?

\bibliographystyle{alpha}
\bibliography{specialization}

\end{document}